\documentclass{amsproc}
\usepackage{amsmath}
\usepackage{amssymb}
\usepackage{amsthm}
\usepackage{amscd}
\usepackage{color}
\usepackage[mathscr]{eucal}
\usepackage{url}
\usepackage{ascmac}
\usepackage{array}
\usepackage[all]{xy}
\usepackage{multirow,eepic}
\usepackage{extarrows}
\usepackage[T1]{fontenc}
\usepackage{lmodern}
\usepackage{bm}
\usepackage{cancel}
\DeclareMathAlphabet{\mathsfit}{\encodingdefault}{\sfdefault}{m}{sl}
\SetMathAlphabet{\mathsfit}{bold}{\encodingdefault}{\sfdefault}{bx}{sl}

\newcommand{\MG}{SL_2(\Z)}

\newcommand{\C}{\mathbb{C}}
\newcommand{\Z}{\mathbb{Z}}

\newcommand{\R}{\mathbb{R}}
\newcommand{\N}{\mathbb{N}}
\newcommand{\Q}{\mathbb{Q}}

\renewcommand{\H}{\mathbb{H}}

\newcommand{\unit}{\mathbf{1}}

\newcommand{\frakg}{\mathfrak{g}}

\newcommand{\ch}{\operatorname{ch}}
\newcommand{\sch}{\mathfrak{ch}}

\newcommand{\wt}{{\rm wt}}
\renewcommand{\=}{\,=\,}
\renewcommand{\H}{\mathbb{H}}

\newcommand{\NO}{\,{\raise0.25em\hbox{$\mathop{\hphantom{\cdot}}%
\limits^{_{\circ}}_{^{\circ}}$}}\,}% The definition of the normal ordering.

\newcommand{\vsh}{\vec{\mathsfit{h}}}
\newcommand{\vh}{\vec{h}}

\newcommand{\SL}{SL_2(\R)}

\makeatletter
\renewcommand{\@biblabel}[1]{#1.}
\makeatother

\newcommand{\G}{\Gamma}

\newcommand{\vf}{\mathbb{F}}

\newcommand{\sd}{\mathfrak{d}}
\makeatletter
\renewcommand{\@biblabel}[1]{#1.}
\makeatother
\newcommand{\mn}{\medskip\noindent}
\newcommand{\bn}{\bigskip\noindent}
\newcommand{\sn}{\smallskip\noindent}
\newcommand{\ms}{\medskip}
\newcommand{\bt}{\begin{theorem}}
\newcommand{\et}{\end{theorem}}

\newcommand{\DEF}{\textbf{Definition.}\;}

\def\fg{\mathfrak{g}}

\makeatletter
\newcommand{\textarrow}[2][1]
  { \settowidth{\@tempdima}{#2}
    \stackrel{#2}
             {\makebox[#1\@tempdima][l]{\rightarrowfill}}
  }
\makeatother

\newcommand{\Rem}{\textbf{Remark.}\;}
\newcommand{\Rems}{\textbf{Remarks.}\;}

\def\be{\begin{equation}}   \def\ee{\end{equation}}     \def\bes{\begin{equation*}}    \def\ees{\end{equation*}}
\def\ba{\be\begin{aligned}} \def\ea{\end{aligned}\ee}   \def\bas{\bes\begin{aligned}}  \def\eas{\end{aligned}\ees}

\renewcommand{\=}{\,=\,}

\newcommand{\AL}[3]{L_{#1_{#2},\,#3}}

\newcommand{\vecr}{\vec{r}}
\newcommand{\vech}{\vec{h}}
\newcommand{\vr}{\vec{r}}
\newcommand{\Gauss}[2]{{}_{#1}F_{#2}}

\newtheorem{theorem}{Theorem}[section]
\newtheorem{proposition}[theorem]{Propositiom}
\newtheorem{lemma}[theorem]{Lemma}

\theoremstyle{definition}

\theoremstyle{remark}

\numberwithin{equation}{section}

\begin{document}
%12.12.2018. 

\title{Vertex operator algebras with central charge~8 and 16}

\author{Kiyokazu Nagatomo }
\address{Department of Pure and Applied Mathematics\\ 
Graduate School of Information Science and Technology\\
Osaka University, Suita, Osaka 565-0871, JAPAN}
%\curraddr{Department of Mathematics and Statistics,
%Case Western Reserve University, Cleveland, Ohio 43403}
\email{nagatomo@math.sci.osaka-u.ac.jp}
\thanks{The first author was supported in part by JSPS KAKENHI Grant Number 17K04171, 
International Center of Theoretical Physics, Italy, and Max Planck institute for Mathematics, Germany.}

\author{Geoffey Mason}
\address{Department of Mathematics\\
 University of California\\
 Santa Cruz, CA 95064, USA}
\email{gem@cats.ucsc.edu}
\thanks{The second author is partially supported by Simons Foundation \#427007.}
%Department of Mathematics, University of California, Santa Cruz, CA 95064, USA
%e-mail: gem@cats.ucsc.edu

\author{Yuichi Sakai}
\address{Multiple Zeta Research Center\\
Kyushu University\\
744, Motooka, Nishi-ku, 
Fukuoka, 819-0395, JAPAN}
\email{dynamixaxs@gmail.com}
\thanks{The third author was supported in part by JSPS KAKENHI Grant numbers 18K03215 and 16H06336.}

%    General info
\date{January 1, 1994 and, in revised form, June 22, 1994.}



\begin{abstract}
In this paper we will partly classify spaces of characters of vertex operator algebras with central charges 8 and 16, 
whose spaces of characters are 3-dimensional and each space of characters forms a~basis of the space of solutions of
a~third order monic modular linear differential equation with rational indicial roots.  

Under the condition that the dimension of the lowest weight space, say $d$,  divides the dimension of the space of weight $d+1$, 
we clarify that the spaces of characters of those vertex operator algebras coincide with the spaces of characters 
of lattice vertex operators associated with integral lattices $\sqrt{2}E_8$ or the affine vertex operator algebra of type $D_{20}^{(1)}$ 
(denoted by $L_{D_{16}^{(1)},1}$
with level 1 for $c=8$, and the Barnes--Wall lattice $\Lambda_{16}$, the affine vertex operator algebras of type $D_{16}^{(1)}$ with level 1 and type $D_{28}^{(1)}$ with level 1 for $c=16$. (The central charge of the affine vertex operator algebras of type 
$D_{28}^{(1)}$ with level 1 is 28, but the space of characters satisfies the differential equations for $c=16$.)

Moreover,  if we suppose a mild condition on characters of vertex operator algebras, then it uniquely determines 
(up to some isomorphisms) the spaces of characters of the lattice $\sqrt{2}E_8$ and the Barnes--Wall lattice $\Lambda_{16}$, respectively.

Finally, the reason why vertex operator algebras with central charges 8 and 16 are intensively studied is that there are solutions such that they does not depends on extra parameters (which represent conformal weights). More precisely, $E_4/\eta^8$ and $E_4^2/\eta^{16}$ are solutions of MLDEs for central charges 8 and 16, respectively. This fact is well understood 
by using hypergeometric function ${}_3F_2$. Hence we cannot apply our standard method to classify vertex operator algebras in which we are interested.

In appendix we classify $c=4$ vertex operator algebras with the same conditions mentioned above.
\end{abstract}

\maketitle

\section*{Introduction}
The principal goal of this paper is to list sets of triplets of Fourier expansions which coincide with the sets characters 
of vertex operator algebras (simply VOA) with central charges 8 and~16. The motivation for doing this was that the lattice VOAs 
associated with $\sqrt{2}E_8$ $(c=8)$ and the Barnes--Wall lattice $(c=16)$ have three simple modules 
and the spaces of characters are 3-dimensional. However, we found that VOAs with central charges 8 and 16 are very interesting 
by the reason stated in the following paragraph. It would be a very important problem if we could classify all lattice VOAs. 
However, this would be difficult since, first of all, even lattices with non-degenerate symmetric forms are not classified 
and secondly, non-isomorphic lattice could define isomorphic lattice VOAs. Therefore we introduce a new concept of ``isomorphisms'' 
(pseudo-isomorphism, weak pseudo-isomorphism) and try to classify VOAs whose sets of characters have central charges $c=8$ 
and $c=16$ under these ``isomorphisms.''

A general form of modular linear differential equation (MLDE) of third order has two parameters $c$ and $h$, 
which correspond to central charges and conformal weights if the MLDE has characters of a VOA as a solution. 
A surprising aspect of the solutions of MLDEs of third order with central charges 8 and 16 is that there are solutions 
without depending on $h$, which are $E_(q^2)/\eta(q)^8$ and $E_4(q)^2/\eta(q)^{16}$, respectively. 
(In fact, $E_(q^2)/\eta(q)^8$ is a~modular function of weight~0 on~$\Gamma(6)$ since $E_4(q^2)$ and $\eta(q)^8$ are modular forms
of weight 4 on $\Gamma_0(2)$ and $\Gamma(3)$, respectively, while $E_4(q)^2/\eta(q)^{16}$ is a modular function of weight 0 
on $\Gamma(3)$ since $E_4(q)^2$ and $\eta(q)^{16}$ are modular forms of weight 8 on $\MG$ and $\Gamma(3)$, respectively.) are solutions of MLDEs for any $h$ and $c=8$ and $6$, respectively. This fact can been seen well from a point of view of writing solutions
in terms of the hypergeometric series $\Gauss{3}{2}$. (See~\eqref{fnhi} in Section~\ref{HGF}.) One can show that 
\begin{equation*}
K^{-c/24}\Gauss{3}{2}\bigl(-c/24,\,h-c/24+1/3,\,c/12+1/2-h+2/3\,;\,1-h,1/2 - c/8 + h\,;\,K\bigr)
\end{equation*}
is independent of $h$ if $c=8$ and 16, where $K=1738/j$ and $j$ is the $j$-function.

A basic tool for investing our problems is the concept of \textit{vector-valued modular forms} (simply VVMF) 
which was introduced by Mason in~\cite{M}. One of advantages of VVMFs is that under some linear condition on leading powers 
(which we call \textbf{exponents}) of Fourier expansions of component functions of a VVMF forms a basis of the space of solutions 
of MLDEs whose singularities are located at the origin $q=e^{2\pi i\tau}=0$ ($\tau=+\infty\sqrt{-1}\pi$) and are regular. 
Similar differential equations were previously studied by Milas in~\cite{AM1}--\cite{AM4},  Adamovic and Milas in~\cite{DM}.

Let $\rho$ be a representation of the full modular group $\rho:\G_1=\MG\rightarrow GL_n(\C)$ and 
$\vf={}^t(f_1,\,f_2,\,\ldots\, f_n\,)$ a column vector-valued holomorphic function on the Poincar\'{e} upper-half plane $\H$
(with Fourier expansions) equipped with a property 
\begin{equation*}
\vf|_k\gamma\={}^t(f_1|_k\gamma,\,f_2|_k\gamma,\,\ldots,\, f_n|_k\gamma)\=\rho(\gamma)\vf
\end{equation*} 
for every $\gamma\in\G_1$, which is called a \textit{vector-valued modular form} (VVMF). It is known by Zhu~\cite{Z}
that the space of characters of a simple, rational VOA of \textit{finite type} (which is usually called $C_2$-cofinite) is invariant 
under the slash zero action $|_0$ of weight zero on the group $\G_1$. Therefore each basis of the space of characters gives a VVMF 
of weight zero on~$\G_1$. If it forms a basis of a linear differential equations with regular singularity only at the origin, 
we can apply the ``generalized'' Frobenius method to determine each character {\it uniquely up to a~scalar factor} 
by its conformal weight as long as incidial roots does not have integral differences. 
However, there are examples that such differential equation does not exist, for instance, 
the affine VOA associated with the Lie algebra $A_2$ with level 3 (see~\cite{AKNS} for instance). 
It is noteworthy that even if indicial roots have integral differences, the differential equations does not have \textit{logarithmic solutions}, i.e., solutions depending on not only $q^{2\pi i\tau}$ but $\tau$ because $V$ is rational (strictly speaking since Zhu's algebras are simple)
as known by Zhu~\cite{Z}. Therefore the spaces of characters of VOAs provide very special examples of the spaces of solutions 
of MLDEs in the sense that logarithmic solutions  do not appear even if its indicial roots have integral differences. 
  
Let $\sd(f)=q\dfrac{d}{dq}(f)-\dfrac{k}{12}E_2(f)$ be the \textit{Serre derivative}, where $f$ is a meromorphic function on $\H$
and $E_2$ is the normalized (quasi) Eisenstein series (of weight two). If $f$ is a modular form of weight $k$, then $\sd(f)$ 
is a modular form of weight $k+2$. Then a linear differential equation of $n$th order 
\begin{equation*}
\sd^n(f)+\sum_{j=2}^nP_{2j}\sd^{n-j}(f)\=0
\end{equation*}
is called a \textbf{modular linear differential equation} (simply MLDE), if $P_{2j}$ is a holomorphic modular forms of weight $2j$ 
on $\G_1$. The coefficient of $\sd^{n-1}$ is~0 since there does not exist a holomorphic modular form of weight two on $\G_1$. 
Let $\vf={}^t\,(f_1,\,f_2,\,\dotsc,\,f_n)$ be a VVMF with weight $k$, whose component functions have Fourier expansions $f_i=c_i^0q^{\lambda_i}+c_i^1q^{\lambda_i+1}+\cdots$ $(c_i^{(0)}\neq0,\,\lambda_i\in\R)$ are linearly independent. 
Then, without loss of generality we may assume that $\lambda_1>\lambda_2>\cdots>\lambda_n$.
(This procedure is said to be the \textbf{normalization}.) Suppose that $n(n+k-1)-12\sum\lambda_i=0$ 
which we call the \textit{non-zero Wronskian condition}. (We introduced this terminology since this equality is equivalent to 
the non-vanishing property of the \textit{modular Wronskian} of the component functions of a~VVMF.) 
Then the set $\{\,f_1,\,f_2,\,\dotsc,\,f_n\,\}$ forms a basis of the space of solutions of a~MLDE (see~\cite{CM} and~\cite{M} for proofs).

One of main purposes of this paper is to show that MLDEs play an important rule when they are applied to the theory of VOAs 
so far to study the spaces of characters. It is clear that the classification of VOAs is the one of the important problems 
of VOAs, however, which would not be easy.  We propose instead the \textit{classification of spaces of pseudo-characters}.  
From the point of view of \textit{conformal field theory} (CFT), the main problems is to determine the system 
of current correlation functions. We therefore try to classify the space of characters, though it might be far from the classification 
of CFTs. Nevertheless, solving this problem contributes the classification of VOAs since there are many cases that if the sets of characters are same, then the corresponding VOAs are isomorphic. (See the minimal models~\cite{ANS} and \cite{ANS2} for examples.)

Let $V$ be a rational, simple VOA of CFT and finite type with a central charge~$c$. Then the space $\sch_V$ of characters 
of simple modules of $V$ is invariant under the usual slash zero action of weight zero on $\MG$. We will show under a mild condition
(\textit{we do not impose neither rationality nor finite type on $V$}\,) that the $\sch_V$ coincides with the space of characters 
of the lattice VOA associated with the integral even lattice $\sqrt{2}E_8$ or the affine VOA of type $D_{20}$ with level one 
if it has the central charge $c=8$, and the \textit{Barnes--Wall lattice} $\Lambda_{16}$ (cf.~\cite[\S4.10]{CS} and~\cite{G}),  
the affine VOAs of type $D_{16}$ with level one and type\footnote{
The central charge of the affine VOAs of type $D_{28}$ with level one is 28. However, its space of characters of simple modules 
is 3-dimensional and coincides with the one of $\Lambda_{16}$.
} 
$D_{28}$ with level one for $c=16$. 

Finally, we impose the extra condition on characters that the second coefficients of characters of vacuum modules coincide with 
ranks of $\sqrt{2}E_8(\,=8)$ and Barnes--Wall lattice $\Lambda_{16}(\,=16)$. Then the both spaces of characters of lattice VOAs are uniquely determined, however, this condition might be slightly strong.

In Section~\ref{thirdMLDE} we review the concepts of VVMFs and MLDEs. The notion of ``the non-zero Wronskian condition'' 
is also defined here. The group of symmetries of MLDEs of third order is discussed in Section~\ref{S.SYM}. 
Section~\ref{S.Lattice} recalls the lattice VOAs associated with the root lattice~$D_n$, which will be used later.
The VOAs with central charge 8 and 16 are studied in Section~\ref{S.cc=8} and Section~\ref{S.cc=16}, respectively,  
and the main theorems of this paper are stated and proved in the these sections. Moreover, in Appendix we give 
the precise expression of the Fourier coefficient mentioned in the footnote in Section~\ref{S.cc=16}, and will study 
$c=4$ case since $c=4n$ ($n\in\N$) cases seems to be interesting.

\mn
\textbf{Warning.} (a) We have often used the terminologies ``central charge'' and ``conformal weight.''
But it is an abuse of terminologies because $c$ and $h$ are a central charges and a conformal weight
if and only if there exists a VOA whose central charge is $c$ and conformal weight are $h$.
Therefore we should call, for instance,~$c$ and $h$ a ``formal'' central charge and a ``formal'' conformal weight, respectively.
However, in many cases, it is not difficult to distinguish formal one and not formal one. Hence we call them just central charges
and conformal weights, which will not cause any confusions.\\
(b) The paper does not intend the classification of vertex operator algebras by isomorphisms, but the one of them by the concept of 
(weakly) pseudo-isomorphisms, which is defined at very end of Section~\ref{S.Lattice}

\section{Modular linear differential equations of third order}\label{thirdMLDE}
Let $\vf={}^t\,(f_1,\,f_2,\,\dotsc,\,f_n)$ be a vector-valued holomorphic (on $\H$) modular form of weight~$k$, 
whose component functions have Fourier expansions, i.e. $f_i=c_iq^{\lambda_i}+O(q^{\lambda_i+1})$ ($c_i\neq0$). 
We suppose that the functions $f_1,\,f_2,\,\dotsc,\, f_n$ are linearly independent and indices $\lambda_i$ are \textit{rational numbers} such that\footnote{
By the successive fundamental column transformations we can always arrange indices as this order.
}
$\lambda_1>\lambda_2>\cdots>\lambda_n$, where $\lambda_i$ is the lowest power of Fourier expansions of $f_i$ of $\vf$. 
We now suppose that $\lambda_i$ satisfies $n(n+k-1)-12\sum_{i=1}^n\lambda_i=0$. 
Then the set $\{f_i\}$ forms a basis of the space of solutions of a MLDE as being proved by Marks and Mason 
(see~\cite[Theorem~2.8]{CM} and \cite[pp.~384]{M} for proofs). Therefore we say that such VVMFs satisfy 
the \textbf{non-zero Wronskian condition}\footnote{
This concept derives from that the Wronskians of the component functions of VVMFs are not zero everywhere.
}. 
The general form of MLDEs of third order is written as (using a derivation $'$ defined below)
\begin{equation}\label{E.Third.EQ1}
f'''-\frac{1}{2}E_2f''+\Bigl(\frac{1}{2}E'_2-xE_4\Bigr)f'+yE_6f\=0\,,\qquad '\=\frac{1}{2\pi i}\frac{d}{d\tau}\,,
\end{equation}
where $E_k$ is the normalized Eisenstein series of weight $k$ on the group $\G_1=\MG$ and $x,\,y$ are real numbers.
(The parameters $x$ and $y$ can be complex numbers, however, in this paper it must be rational numbers 
since we are interested in rational VOAs, whose central charges and conformal weights are rational.)

Let $V$ be a VOA with a rational central charge $c$, whose space of characters is 3-dimensional, and let $\{0,\,h_1,\,h_2\}$ 
be the set of conformal weights whose elements are rational numbers. Suppose that the set $\{-c/24,\,h_1-c/24,\,h_2-c/24\}$ 
of rational numbers consists of \textbf{indicial roots} (the roots of the indicial equation) of the MLDE associated with 
a VVMF satisfying the non-zero Wronskian condition. Then the set of conformal weights obliviously can be written as 
$\{0,\,h,\,c/8-h+1/2\}$ with a rational number $h$. (This expression is, in fact, not unique, we have two choices either 
$h=h_1$ or $h_2$). Therefore we one can confirm that every MLDE of our interest is written as
\begin{multline}\label{E.Third.EQ2}
f'''-\frac{1}{2}E_2f''+\Bigl\{\frac{1}{2}E'_2-\Bigl(h^2-\frac{h}{2}-\frac{ch}{8}+\frac{c^2}{192}+\frac{c}{24}\Bigr)E_4\Bigr\}f'\\
+\frac{c}{24}\Bigl(\frac{c}{12}+\frac{1}{2}-h\Bigr)\Bigl(h-\frac{c}{24}\Bigr)E_6f\=0\,.
\end{multline}
It is clear by the very definition of the indicial equation that the set of indices of the MLDE~\eqref{E.Third.EQ2} 
is~$\left\{-c/24,\,h-c/24,\,c/12-h+1/2\right\}$.

Suppose that~\eqref{E.Third.EQ2} has a solution of the form $f=q^\lambda\bigl(\sum_{n=0}^{\infty}a_nq^{n}\bigr)$ with $a_0=1$ 
and $\lambda$ is one of the indicial roots $-c/24$, $h-c/24$ and $c/12-h+1/2$. Then using the Frobenius method, one can determine Fourier coefficients $a_n$ $(n\in\N)$ uniquely by the recursive relation
\begin{equation}\label{thirdMLDE.REC}
\begin{split}
&\left(n+\lambda+\frac{c}{24}\right)\left(n+\lambda+\frac{c}{24}-h\right)
\Bigl(n+\lambda-\frac{c}{12}-\frac{1}{2}+h\Bigr)a_n\\
&\=\sum_{i=1}^{n}\biggl\{(\lambda-i+n)\Bigl(12(2i-\lambda-n)\sigma_1(i)\\
&\phantom{\qquad\times\sum_{i=1}^{n}}+\frac{5}{4}\left(c^2+8c-24hc-96h+192h^2\right)
\sigma_3(i)\Bigr)\\
&\phantom{\qquad\times\sum_{i=1}^{n}\biggl\{(\lambda-i+n)}-\frac{7}{96}c(c-24h)(c-12h+6)\sigma_5(i) \biggr\}a_{n-i}
\end{split}
\end{equation}
as far as the factors of the left-hand side is not zero, i.e., the indicial roots do not have integral difference, where the $\sigma_m(i)$
is the sum of the $m$th powers of the divisors of~$i$. The denominator (the left-hand side of the recursions) of $a_n$ $(n\in\N)$ 
are zero if and only if the following situations simultaneously happen (see Table~\ref{Table_Indefinite}).
\begin{table}[htbp]
\begin{center}
\begingroup
\renewcommand{\arraystretch}{1.2}
\begin{tabular}{c|ll}
$\lambda$ & The possibility of indefinite values of $\{a_n\}_{n\in\N}$ \\ 
\hline
$-c/24$&$h\in\N,\, c/8-h+1/2\in\N$\\ 
$h-c/24$&	$-h\in\N,\,c/8-2h+1/2\in\N$\\
$c/12+1/2-h$&$h-c/8-1/2\in\N,\,2h-c/8-1/2\in\N$	
\end{tabular}
\endgroup

\mn
\caption{Indices and the possibility of indefinite values of the coefficients}
\label{Table_Indefinite}
\end{center}
\end{table}

\mn
\Rems
(a)
If $c=8$ and  $\lambda=-1/3$, then denominators appeared in $a_1$ and $a_2$ are zero if $h=1$ and $h=1/2$, respectively. 
However, it follows that $a_1$ and $a_2$ are 248 and 4124, respectively,  after canceling common factors as 
\begin{align}\label{E.Sspecial1}
&a_1\=\frac{992 \cancel{(h-1)(1-2h)}}{4\cancel{(h-1)(1-2h)}}\=248\,,\\
&a_2\=\frac{131968 
\cancel{(h-2)(h-1)}\cancel{(2h-1)(1+2h)}}{32\cancel{(h-2)(h-1)}\cancel{(2h-1)(1+2h)}}\=4124\,.\label{E.Sspecial1}
\end{align}
(b)
For $c=8$ and $\lambda=h-1/3$ and $7/6-h$ the $a_n$ never become indefinite if $-h,\,3/2 - 2 h\not\in\N$ 
and $h-3/2,\, 2h-3/2\not\in\N$ by Table~\ref{Table_Indefinite}.
\\
(c)~If $\lambda=-c/24$, the coefficient $a_1$~is well defined when~$h\neq0,\,1,\, c/8\pm1/2$. 
For $c=8$~and~16, the $q$-series which we obtain are 
\begin{align}
&E_4/\eta^8\=q^{-1/3}\bigl(1 + 248 q + 4124 q^2 + 34752 q^3 + 213126 q^4+\cdots\bigr)\,,\label{E.c8E4}\\
&E_4^2/\eta^{16}\=q^{-2/3}\big(1 + 496 q + 69752 q^2 + 2115008 q^3 + 34670620 q^4+\cdots\bigr)\,,\label{E.c16E4}
\end{align}
respectively. In fact, the functions~\eqref{E.c8E4} and~\eqref{E.c16E4} satisfy~\eqref{E.Third.EQ2} for $c=8$ and $16$ 
for any~$h\in\Q$, respectively,  which can be verified by using Ramanujan's identities 
\begin{equation*}
24\eta'\=E_2\eta\,, \quad 12E_2'\=E_2^2-E_4\,, \quad 3E_4'\=E_2E_4-E_6\,,\quad 2E_6'^\=E_2E_6-E_4^2\,,
\end{equation*}
where $\eta(q)=q^{1/24}\prod_{n=1}^{\infty}(1-q^n)$ is the Dedekind eta function.  
The MLDE of orde 3 with weight $k$, which has a formal $q$-series solution $f=1+O(q)$, has another expression as 
\begin{equation*}
\sd^3_{k}(f)+\alpha E_4\sd_{k}(f)+\left(\dfrac{k(k+2)(k+4)}{12^3}+\dfrac{k}{12}\alpha\right) E_6 f\=0
\end{equation*}
where $\sd_k(f)=\Bigl(f'-\dfrac{k}{12}E_2\Bigr)f$ is the Ramanujan--Serre derivtive, the iterated derivation is defined by 
$\sd_{k}^{i}(f)=\sd^{i-1}_{k+2(i-1)}\circ\sd_{k+2(i-2)}(f)$ $(i\geq1)$, $\sd^{0}_{*}(f)=f$ and $\alpha\in\C$. 

Then for $k=4$ and~$8$, $E_4$ and $E_4^2$ are solutions of this MLDE of weight 4 and~8, respectively, for any complex number 
$\alpha$. (This fact is explained in William Stein's free online textbook "Computing with Modular Forms.")
Therefore it may be notable that there are solutions of~\eqref{E.Third.EQ2} \textit{do not depend on the variable~$h$}, 
whose Fourier expansions with indices~$-1/3$ ($c=8$) and~$-2/3$ ($c=16$), respectively.

\mn
\Rem
We say that any solution of a MLDE with non-negative integral coefficients $a_n$ is \textbf{of character type}, particularly, 
is \textbf{of vacuum character type} if $a_0=1$. In~\cite{NS} we have used the fact that the solution
$f_1=q^{-c/24}\left(\sum_{n=0}^{\infty}a_nq^n\right)$ with $a_0=1$, which is of vacuum character type 
of $V$ really does depend on the parameter~$h$, more precisely, $a_1$ is a non-trivial rational function in $c$ and~$h$, 
and we were able to obtain values of $h$ which give all solutions of character type of~\eqref{E.Third.EQ2}. 
The fact that the coefficients of $\ch_V$ are non-negative integers played an important role to determine~$h$. 
\textit{However, this approach does not work on the present case since\footnote{
That is the solution with the index $-c/24$ obtaining by the Frobenius method.
} $f_1$ does not depend on the variable~$h$.} Therefore we will take an altanative approach in Section~\ref{S.cc=8} and~Section~\ref{S.cc=16}. But, before describing the results, we introduce symmetries of~\eqref{E.Third.EQ2}, 
which plays a crucial rule to restrict the possible values of conformal weights~$h$.

\section{Symmetries of modular linear differential equations of third order}\label{S.SYM}
It is shown easily that the MLDE~\eqref{E.Third.EQ2} does not change under the action of the group
(which is isomorphic to $S_3$) generated by two affine isomorphisms $\lambda$ and $\mu:\C^2\rightarrow\C^2$ defined by
\begin{equation*}
\lambda(c,\,h)\=(c-24h,\, c/8-2h+1/2)\quad\text{and}\quad\mu(c,\,h)\= (c-24h,\, -h)\,,
\end{equation*}
which have order 3 and~2, respectively. Suppose that~\eqref{E.Third.EQ2} has solutions of the form 
\begin{align*}
&f_1\=q^{-c/24}+O(q^{1-c/24})\,,\\
&f_2\=q^{h-c/24}+O(q^{1+h-c/24})\,,\\
&f_3\=q^{c/12+1/2-h}+O(q^{c/12+3/2-h})\,,
\end{align*}
which are obtained by the \textit{Frobenius method}, i.e., the coefficients $a_n$ are determined recursively 
by~\eqref{thirdMLDE.REC}. 

\mn
\textbf{Warning.} 
Throughout this paper, we use $f_1$, $f_2$ and $f_3$ to denote solutions of~\eqref{E.Third.EQ2} if they exist, 
whose indices $-c/24$,  $h-c/24$ and $c/12+1/2-h$, respectively, i.e, we always work under the condition given 
in Table~\ref{Table_Indefinite}.

\ms
Then the action of the group $S_3$ defined above gives permutations of the set of solutions gives
\begin{equation*}
\lambda:(f_1,\, f_2,\, f_3)\longmapsto(f_2,\, f_3,\, f_1)\,,\quad \mu:(f_1,\, f_2,\, f_3)\longmapsto(f_1,\, f_3,\, f_2)\,.
\end{equation*}
Therefore for any two numbers $c$ $h$ there exist $\sigma,\,\tau\in S_3$ such that $\sigma(c)$ and $\sigma(h)$ can be considered 
as a~central charge and a conformal weight, respectively, (where we supposed that $c$ and $h$ are originally central charges 
and conformal weights, respectively) and there are three ways that we can think the elements 
of $\left\{-c/24,\,h-c/24,\,c/12-h+1/2\right\}$ as a central charge and conformal weights as being listed in Table~\ref{SYS.TL.CW}.
\begin{table}[htbp]
\begin{center}
\begingroup
\renewcommand{\arraystretch}{1.2}

\begin{tabular}{l|l}
Central charge & Conformal weights \\ \hline
$c$&  $0$, $h$, $c/8+1/2-h$\\ 
$c-24h$ & $-h$, $0$, $c/8+1/2-2h$\\ 
$24h-2c-12$ &$h-c/8-1/2$, $2h-c/8-1/2$, $0$
\end{tabular}
\endgroup

\mn
\caption{Central charge and conformal weights}\label{SYS.TL.CW}
\end{center}
\end{table}

\noindent
For instance, the space of solutions of~\eqref{E.Third.EQ2} with the central charge $c$ and the conformal weight $h$ 
is same as the one with the central charge $c-24h$ and the conformal weight 0 (see the second row of the Table~\ref{SYS.TL.CW}).

Suppose that $(h-1)(c-8h-4)\neq0$. Then we have $\neq1$, otherwise, since $\lambda^2(c,h)=(-2c+12, (-c+4)/8)\neq(c,1)$.
Therefore $h=1$ and then $(h-1)(c-8h-4)=0$ which is a contradiction.Then a short calculation shows that the second coefficient
$m$ of $f_1$ is given by
\begin{equation*}
m\=\frac{-c^3+31c^2h-7c^2-248ch^2+124ch-4 c}{(h-1)(c-8 h-4)}\,
\end{equation*}
which is rewritten as 
\begin{equation}\label{SYM.Po}
8(m-31c)h^2-(4 + c)(m-31c)h-c^3-7c^2+(m-4)c-4m\=0\,.
\end{equation}
The above equation has two solutions 
\begin{equation}\label{SYM.H}
h\= \frac{1}{16}\left(c+4\pm\frac{\sqrt{\mathsf{D}}}{m-31c}\right)\,
\end{equation}
if and only if $m\neq31c$. If $m=c/31$, then we have $m=0,\,248,\, 496$.

Suppose that $c$ is an \textit{integer}. Then $h$ is a rational number if and only if the discriminant
\begin{equation}\label{SYM.D}
\begin{split}
\mathsf{D}&\=(m-31c)\bigl(m(c-12)^2+c(c^2-24c-368)\bigr)\\
&\=(c-12)^2\left(m-\frac{c(15c^2-360c+2416)}{(c-12)^2}\right)^2-\frac{16^2c^2(c-8)^2(c-16)^2}{(c-12)^2}
\end{split}
\end{equation}
is the square $d^{\,2}$ of an~integer~$d$. (Moreover,  if~$c$ is a~rational number, $\mathsf{D}$ is the~square 
of a~rational number.) Then~\eqref{SYM.D} is equivalent to (by using $d$\,)
\begin{equation*}
\begin{split}
&16^2c^2(c-8)^2(c-16)^2\\
&\qquad\qquad\=\left((c-12)\left((c-12)m-d\right)-c(15c^2-360c+2416)\right)\\
&\qquad\qquad\qquad\qquad\qquad\qquad\times\left((c-12)\left((c-12)m+d\right)-c(15c^2-360c+2416)\right)\,.
\end{split}
\end{equation*}
Since $c$ is an integer, one can verify that
\begin{equation*}
(c-12)\left((c-12)m\pm d\right)-c(15c^2-360c+2416)
\end{equation*}
are divisors of $16^2c^2(c-8)^2(c-16)^2$. Then each $c$ defines an integer $m=c/31$ and then rational numbers $h$ by~\eqref{SYM.H}.

If $\sqrt{\mathsf{D}}$ is an~integer, obviously there are two possible rational values of~$h$ given by~\eqref{SYM.H}, 
which we denote by $h_1$ and~$h_2$ (or $h_2$ and~$h_1$), respectively. Now, say $h_1=(c+4+\sqrt{\mathsf{D}}/(m-31c))/16$. 
Then it follows that  
\begin{equation*}
\mu(h_1)\=\frac{c}{8}+\frac{1}{2}-h_1\= \frac{c+4-\sqrt{\mathsf{D}}/(m-31c)}{16}\=h_2\,.
\end{equation*}
Therefore we have the following 
\begin{lemma}\label{SYM.INV}
The set of values of $h$ in {\rm Table~\ref{SYS.TL.CW}} is closed under the action of $\mu\in S_3$, where $\mu(c,\,h)\= (c-24h,\, -h)$.
\end{lemma}

We use the following useful 

\mn
\textbf{Convention.}
Two rational numbers $c$ and $h$ are considered as a central charge and a conformal weight, respectively, whenever we write 
$(c,\,h)\in\Q^2$. The image of $(c,\,h)$ of a symmetry in $S_3$ is denoted by Sans-Serif letters as $(\mathsfit{c},\,\mathsfit{h}\,)$.

\mn
\DEF
Let $c$, $\mathsfit{c}$ and $h$, $\mathsfit{h}$ be complex numbers. We say that the ordered pairs of complex numbers 
$(c, \,h)$ and $(\mathsfit{c},\,\mathsfit{h}\,)$ are \textbf{equivalent} if there is an element $\sigma$ 
of $\langle\lambda,\mu\rangle\,(\simeq S_3)$ such that $\sigma(c,h)=(\mathsfit{c},\, \mathsfit{h})$, 
We write the sets of corresponding conformal weights by $\vh=(h_1, \,h_2)$ and $\vsh=(\mathsfit{h}_1, \mathsfit{h}_2)$, respectively.

\mn
\Rems
(a) Since $\lambda$ and $\mu$ are isomorphisms the equivalent relation $\thicksim$ is well defined.\\
(b) For $c=8$ or $16$, we can suppose that $(h-1)(c-8h-4)\neq0$.

\section{Theta series associated to the lattice $D_n$}\label{S.Lattice}
In this section (for the reader's convenience) we give the explicit forms of characters of simple modules over the lattice VOA associated with the non-degenerate even integral lattice $L=D_n$ $(n\geq4)$ in terms of the lattice theta series, which will be used later. 

Let $e_i=(\delta_{ij})\in\Z^n$ for $1\leq i\leq n$, where $\delta_{ij}$ denotes the Kronecker delta. 
We define the lattice $L=D_n$ by 
\begin{equation*}
L\=\bigl\{x=(x_1,\,x_2,\,\dotsc\,,x_n)\in\Z^n\;\Bigl|\;\sum_{i=1}^{n}x_i\equiv0\,\bmod{2}\bigr\}\,,
\end{equation*}
which is also written as $L=\langle\, e_1-e_2,\, e_2-e_3,\,\ldots,\,e_{n-1}-e_{n},\, e_{n-1}+e_{n}\,\rangle_{\Z}$.

The dual lattice of $L$ is defined by $L^\circ=\left\{v\in\Q\otimes_{\Z}L\,\big|\, \langle v,L\rangle\subseteq\Z\;\right\}$, 
where $\langle\;\,,\;\rangle$ is the standard inner product on~$\Q$. Let $\xi=e_1$ and $\eta=(e_1+\dotsb+e_n)/2$. 
Then the quotient vector space $L^\circ/L$ is generated by two residue classes $\xi$ and~$\eta$ 
(cf. see~\cite[pp.\,89, Example~5.3]{HK} for more details). The theta series $\Theta_L$ of the lattice $L$ and the cosets 
of the dual lattice $L^\circ$ are given by 
\begin{align*}
&\Theta_{L}(q)\=\frac{\theta_3(q)^n+\theta_0(q)^n}{2}\=1+2n(n-1)q+\frac{2}{3}n(n^3-6n^2+11n-3)q^2 +\cdots\,,\\
&\Theta_{\xi+L}(q)\=\frac{\theta_3(q)^n-\theta_0(q)^n}{2}\=2nq^{1/2}\Bigl(1+\frac{2}{3}(n-1)(n-2)q+\cdots\Bigr)\,,\\
&\Theta_{\eta+L}(q)\=\frac{\theta_2(q)^n}{2}\=2^{n-1}q^{n/8}\Bigl(1+nq+\frac{1}{2}n(n-1)q^2+\cdots\Bigr)\,,
\end{align*}
respectively, where 
\begin{equation*}
\theta_2(q)\=\sum_{n\in\Z}q^{(n+1/2)^2/2}\,,\quad\theta_3(q)\=\sum_{n\in\Z}q^{n^2/2}\,,\quad\theta_0(q)\=\sum_{n\in\Z}(-1)^nq^{n^2/2}
\end{equation*}
are Jacobi's theta constants which are modular forms of weight $1/2$ on $\G(8)$ (see~\cite[pp.\,118, Section~7.1]{CS} 
for more details). Let $N$ be the maximum divisor of $n$ such that $1\leq N\leq4$, and $N_1$ the number written 
in Table~\ref{THETA_LEVELS}. 
\begin{table}[htbp]
\begin{center}
\begin{tabular}{c|c|c|c|c}
$N\,$&$1$&$2$&$3$&$4$ \\ \hline
$N_1$&$8$&$4$&$8$&$2$
\end{tabular}

\mn
\caption{Levels of powers of Jacobi theta's constants}\label{THETA_LEVELS}
\end{center}
\end{table}

Then, with the help of the theory of modular forms, the functions $\theta_2^n$, $\theta_3^n$ and $\theta_0^n$ $(n>0)$ are modular forms of weight~$n/2$ on~$\Gamma(N_1)$. For the non-expert's convenience the principal congruence subgroup $\G(8)$ 
of level 8 is defined by 
\begin{equation*}
\G(8)\=\left\{
\begin{pmatrix}a&b\\ c&d\end{pmatrix}
\in\G_1\,\bigg|\,
\begin{pmatrix}a&b\\ c&d\end{pmatrix}
\equiv\begin{pmatrix}1&0\\ 0&1\end{pmatrix}\,
\bmod{8}
\right\}\,.
\end{equation*} 
Since characters of simple modules of the VOA $V_L$ are given by $\Theta_{\xi+L}(q)/\eta(q)^{\text{rank}(L)}$ 
and $\text{rank}(L)=n$~(cf.~\cite{FLM}),  one can verify that 
\begin{align}
\begin{split}\label{ch-D1}
&\ch_V(\tau)\=\frac{\Theta_{L}(q)}{\eta(q)^n}\\
&\=\frac{\theta_3(q)^n+\theta_0(q)^n}{2\eta(q)^n}\\
&\= q^{-n/24}\Bigl(1+n(2n-1)q+\frac{1}{6}n(4n^3-12n^2+35n-3)q^2+\cdots\Bigr)\,,
\end{split}\\
\begin{split}\label{ch-D2}
&\ch_{V_{\xi+L}}(\tau)\=\frac{\Theta_{\xi+L}(q)}{\eta(q)^{n}}\=\frac{\theta_3(q)^n-\theta_0(q)^n}{2\eta(q)^{n}}\=2nq^{1/2-n/24}\\
&\qquad\times\Bigl(1+\frac{1}{3}(2n^2-3n+4)q
+\frac{1}{30}(4n^4-20n^3+95n^2-55n+36)q^2+\cdots\Bigr)\,,
\end{split}\\
\begin{split}\label{ch-D3}
&\ch_{V_{\eta+L}}(\tau)\=\frac{\Theta_{\eta+L}(q)}{\eta(q)^n}\\
&\=\frac{\theta_2(q)^n}{2\eta(q)^n}\\
&\=2^{n-1}q^{n/8-n/24}\bigl(1+2nq+n(2n+n)q^2+\cdots\bigr)\,.
\end{split}
\end{align}

The characters~\eqref{ch-D1}, \eqref{ch-D3} and~\eqref{ch-D3} satisfy the monic MLDE of third order
\begin{equation}\label{MLDE_for_D}
f'''-\frac{1}{2}E_2\,f''+\Bigl(\frac{1}{2}E_2'+\frac{n}{48}\Bigl(1-\frac{n}{4}\Bigr)E_4\Bigr)f'+\left(\frac{n}{24}\right)^2
\left(1-\frac{n}{12}\right)E_6\,f\=0\,.
\end{equation}

Let $N_2$ be as in the Table~\ref{ETA_LEVELS}. Since the Dedekind eta function $\eta$ is a modular form of weight $1/2$
on~$\G(24)$, the function $\eta^n$ is of weight $n/2$ on~$\G(N_2)$, where $N$ is the greatest common divisor of $n$ 
and $24$. The function $\eta^n$ $(n=1,\,2,\,3,\,4,\,6,\,8,\,12,\,24)$ is a modular form on $\Gamma(N_2)$ 
$(N_2=24,\,12,\,8,\,6,\,4,\,3,\,2,\,1)$, respectively.

\begin{table}[htbp]
\begin{center}
\begin{tabular}{c|c|c|c|c|c|c|c|c}%\hline
$N$&$1$&$2$&$3$&$4$&$6$&$8$&$12$&$24$ \\ \hline
$N_2$&$24$&$12$&$8$&$6$&$4$&$3$&$2$&$1$%\\ \hline
\end{tabular}
\caption{Levels of powers of the Dedekind eta function ($N=\gcd{(N,\,n})$)}\label{ETA_LEVELS}
\end{center}
\end{table}
Therefore it follows that the characters of affine VOA of type $D_n$~$(n\ge4)$ are modular functions of level $N_3$ 
in Table~\ref{CH_LEVELS}, where $N$ is the greatest common divisor of $n$~and~12. 

\begin{table}[htbp]
\begin{center}
\begin{tabular}{c|c|c|c|c|c|c}%hline
$N$ &$1$&$2$&$3$&$4$&$6$&$12$ \\
\hline
$N_3$&$24$&$12$&$8$&$6$&$4$&$2$%\\ \hline
\end{tabular}
\caption{The levels of characters of affine VOA of type $D_n$ with level~1 $(n\ge4)$ ($N=\gcd{(n,12)}$)}\label{CH_LEVELS}
\end{center}
\end{table}

It is  proved by C.~Dong, X.~Lin and S-H.~NG in~\cite{DLN} (for the reader's convenience $Z_i(v,\tau)$ appeared in the paper 
is a one-point correlation function evaluated at $v\in V$ and $\wt[v]$ is a weight of $V$ but Virasoro element is distinct 
of the one of $V$, in particular, $\chi_i=Z_i(\unit,\tau)$):

\begin{theorem}[{{\cite[Theorem~1]{DLN}}}] \label{T-DLN}
Let $V$ be a rational, self-dual, simple vertex operator algebra of CFT and finite type. Then each $Z_i(v,\tau)$ is a modular form 
of weight $\wt{[v]}$ on a congruence subgroup of $\MG$ of level $n$ which is the smallest positive integer 
such that $n(\lambda_i-c/24)$ is an integer for all  $i$. In particular, each $q$-character $\chi_i$ is a modular function 
on the same congruence subgroup.
\end{theorem}

\mn
\Rem
Since MLDEs which we studied are assumed to come from VOAs, one can suppose that the characters are modular forms 
of weight 0 on the congruence group of level $n$ determined in Theorem~\ref{T-DLN}. However, there are examples of characters 
whose levels of congruence groups are less than~$n$. 

\bn
\DEF 
If the characters of vacuum modules for two VOAs $V$ and $W$ are equal and the set of characters (modules) of two VOAs $V$
and $W$~coincide, we say that these two VOAs are \textbf{pseudo-isomorphic}. Indeed, there is an example that $V$ and $W$ are pseudo-isomorphic but are not isomorphic (see Theorem~\ref{T.c=8.1/2}). Let $\ch_V$ and $\ch_W$ be the spaces of characters 
of $V$ and $W$, respectively. If $\ch_V$ and $\ch_W$ coincide, we say that $V$ and $W$ are \textbf{weakly pseudo-isomorphic}.

\section{Hypergeometric function and modular linear differential equations}\label{HGF}
The main motivation in this section comes from monic MLDEs with $c=8$ and $16$.
We start with a monic MLDE of third order with weight zero, which is expressed by means of the Serre derivative
\begin{equation}\label{MLDE}
(\sd^3+xE_4\sd+yE_6)(f)\=0\,,\quad \sd=\sd_0
\end{equation}
where $x$ and $y$ are complex numbers.

The principal goal of this section is to give an explicit expression of a solution of~\eqref{MLDE} in terms of the hypergeometric 
function ${}_3F_2$, whose coefficients of Fourier expansions are all \textit{non-negative integers}. 
Moreover, we suppose that there is at least a~solution of which leading coefficient is one.

For this purpose we introduce the hypergeometric function ${}_3F_2$ by 
\begin{equation*}
{}_3F_2(a,b,c\,;d,e)\=\sum_{r=0}^\infty\frac{(a)_r(b)_r(c)_r}{(d)_r(e)_r}\cdot\frac{z^r}{r!}\,,
\end{equation*}
where $(x)_r=x(x+1)\cdots(x+r-1)$ for $r>0$ and $(x)_0=1$ for any $x$ is the increasing \textbf{Pochhammer symbol}
(see \cite{FM} for more general hypergeometric series).

Using the hypergeometric function one can define three functions
\begin{equation}\label{HS}
f_i\=K^{r_i}\,{}_3F_2\bigl(r_i,r_i+1/3,r_i+2/3;r_i-r_j+1,r_i-r_k+1;K\bigr)
\end{equation}
for  $i=1,\,2$ and 3, where $K=1728/j$ and $j$ is the $j$-function which has the $q$-expansion $j=q^{-1}+744+\cdots$.

We now let $r_1,\,r_2$ and $r_3$ be exponents of the VVMF $\vf={}^t(f_1,\,f_2,\,f_3)$. %of $\rho(T)$ 
Then it is shown in~\cite{FM} that~\eqref{HS} form a basis of the space of solutions of  the MLDE~\eqref{MLDE} if
$r_i-r_j+1,\, r_i-r_k+1\not\in\Z_{\leq0}$.  The first two terms of Fourier expansion of each $f_i$ are given by
\begin{align}
&f_1/(1728)^{r_1}\=q^{r_1}+\Bigl(1728\frac{r_1(r_1+1/3)(r_1+2/3)}{(1+r_1-r_2)(1+r_1-r_3)}-744r_1\Bigr)q^{r_1+1}+\cdots\,,\\
&f_2/(1728)^{r_2}\=q^{r_2}+\Bigl(1728\frac{r_1(r_2+1/3)(r_2+2/3)}{(1+r_2-r_1)(1+r_2-r_3)}-744r_2\Bigr)q^{r_2+1}+\cdots\,,\\
&f_3/(1728)^{r_3}\=q^{r_3}+\Bigl(1728\frac{r_1(r_3+1/3)(r_3+2/3)}{(1+r_3-r_2)(1+r_3-r_1)}-744r_3\Bigr)q^{r_2+1}+\cdots\,.
\end{align}

The two functions $x$ and $y$ in the MLDE~\eqref{MLDE} are also expressed by using hypergeometric expressions as follows:
\begin{align*}
\begin{split}
x(q)&\=\frac{2E_2(q^2)-E_2(q)}{\eta(q)^4}\\
&\=j(q)^{1/6}{}_3F_{2}\left(-\dfrac{1}{6}\,,\,\dfrac{1}{6}\,,\,\dfrac{1}{2}\,;\,\dfrac{1}{2}\,,\,\dfrac{1}{2}\,;\,\dfrac{12^3}{j(q)}\right)\\
&\=q^{-1/6}(1 + 28 q + 134 q^2 + 568 q^3 + 1809 q^4 +\cdots)\,,
\end{split}\\
y(q)&\=\frac{\eta(q^2)^8}{\eta(q)^8}\\
\begin{split}
&\=j(q)^{-1/3}{}_3F_{2}\Big(\dfrac{1}{3}\,,\,\dfrac{2}{3}\,,\,1\,;\,\dfrac{3}{2}\,,\,1\,;\,\dfrac{12^3}{j(q)}\Big)\\
&\=q^{1/3}(1 + 8 q + 36 q^2 + 128 q^3 + 394 q^4 +\cdots)\,,
\end{split}
\end{align*}
which will be used in Section~\ref{SS.cc=8}.

Setting $r_1=-c/24$, $r_2=h_2-c/24$ and $r_3=h_2-c/24$ one can verify that 
\begin{align}
&f_1/(1728)^{r_1}\=q^{r_1}-c\left(\frac{1728(c-8)(c-16)}{3(1-h_2)(1-h_3)}-31\right)q^{r_1+1}+\cdots\,,\\
\begin{split}
&f_2/(1728)^{r_2}\=q^{r_2}\\
&\qquad\qquad+(h_2-c/24)\left(\frac{1728r_1(r_2+1/3)(r_2+2/3)}{(1+r_2-r_1)(1+r_2-r_3)}-744r_2\right)q^{r_2+1}+\cdots\,,
\end{split}\\
&f_3/(1728)^{r_3}\=q^{r_3}+\left(\frac{1728r_1(r_3+1/3)(r_3+2/3)}{(1+r_3-r_2)(1+r_3-r_1)}-744r_3\right)q^{r_2+1}+\cdots\,.\label{fnhi}
\end{align}
Since $r_1=-c/24$, $r_2=h_1-c/24$ and $r_3=h_1-c/24$, $f_1$, $f_2$ and $f_3$ are rewritten as
\begin{align}
&f_1/(1728)^{r_1}\=q^{r_1}-c\left(\frac{1(c-8)(c-16)}{3(1-h_2)(1-h_3)}-31\right)q^{r_1+1}+\cdots\,,\\
\begin{split}
&f_2/(1728)^{r_2}\=q^{r_2}+(1/2 + c/24 - h)\\
&\qquad\times\left(\frac{c (c-24 h_1-16) (c-24 h_1-8)}{8 (h_1+1)(-h_1+h_2-1)}+31(c-24h_1)\right)q^{r_2+1}+\cdots\,,
\end{split}\\
\begin{split}
&f_3/(1728)^{r_3}\\
&\=q^{r_3}+\left(31 (c-24 h_2)-\frac{c (c-24 h_2-8) (c-8 (3h_2+2))}{8 (h_2+1) (-h_1+h_2+1)}\right)q^{r_2+1}+\cdots\,.\label{fnhi}
\end{split}
\end{align}
Finally, since  $h_1=h$ and $h_2=c/8-h+1/2$, it follows that 
\begin{align}
\begin{split}
&f_1/(1728)^{r_1}\\
&\quad\=q^{r_1}-c\Bigl(\frac{(c-8)(c-16)}{-(1/288) (24 + c - 24 h) (-6 + c - 12 h)}\Bigr)q^{r_1+1}+\cdots\,,\notag
\end{split}\\
\begin{split}
&f_2/(1728)^{r_2} \=q^{r_2}+(1/2 + c/24 - h)\\
&\qquad\times\Bigl(31 (c-24 h)+\frac{c (c-24 h-8) (c-8 (3 h+2))}{(h+1) (c-16 h-4)}\Bigr)q^{r_2+1}+\cdots\,,
\notag
\end{split}\\
\begin{split}
&f_3/(1728)^{r_3}\=q^{r_3}\\
&\qquad+\Bigl(-62 (c-12 h+6)-\frac{32 c (c-12 h+10) (c-12 h+14)}{(c-16 h+12) (c-8 h+12)}\Bigr)q^{r_2+1}+\cdots\,.\label{fnhi}
\end{split}
\end{align}

Now it is clear from the expression of $f_1$  that $c=8$ and $c=16$ are special central charges for the third order MLDE.

\mn
\Rem
Lemma~18 in~\cite{FM1} shows that the following two statements are equivalent.
Let $F$ be a~VVMF. Then 
(1)  At least one of the components of $F$ is a modular form on a congruence subgroup,
(2) All of the components of $F$ are modular forms on a congruence subgroup.

\section{Vertex operator algebras with central charge~$c=8$}\label{S.cc=8}
In this section motivated by the discussions given in section~\ref{HGF} we study VOAs whose central charges are all~8 . 
We now suppose that~\eqref{MLDE} has a solution of the form $f=q^{-1/3}(1+mq+O(q^2))$ (as~$c=8$). 
Substituting $f$ into~\eqref{MLDE} one can find first $x=1/2$ and $y=1/27$ if $m\neq248$. 
However, for $m=248$v we only have the linear relation $18x+54y=5$, and then $(x,\,y)$ is not unique. 
Therefore we first work under the hypothesis $m\neq248$.

\subsection{Generic case}\label{SS.cc=8}
In this subsection we will work  under the hypothesis

\mn
\begin{center}
\textbf{Hypothesis:} $m\neq248$.
\end{center}

\mn
\DEF 
we denote the set of indicial roots and the set of conformal weights by~$\vec{\mathbf{r}}$ and $\vec{\mathbf{h}}$, respectively.

\begin{proposition}[Strum bound]
Let $f$ and $g$ be modular forms of weight $k$ on Fuchsian group $\subset\SL$. Suppose that first $N=kd_\G$ 
$q$-expansions of $f$ and $g$ are equal then we have $f=g$, where $d_\G$ is the index of $\G\subset\SL$

\end{proposition}

Let us denote $f_1$, $f_2$ and $f_3$ the solutions associated with the indicial roots~$\vr=\{\lambda_1,\,\lambda_2,\,\lambda_3\}$.  Then each $\lambda_i$ is a solution of a cubic equation called the \textit{indicial equation} of~\eqref{E.Third.EQ2} is given by 
$t^3- t^2/2 - t/6+1/27=0$ (since $x=1/2$ and~$y=1/27$). 
Hence the set of indicial roots $\vecr=\{-1/3,\, 1/6,\, 2/3\}$ and the set of conformal weights $\vech={0, 1/2, 4/3}$, respectively.

It is obvious that the lattice VOA $V_{\sqrt{2}E_8}$ and the affine VOA $\AL{D}{8}{1}$ are not isomorphic as a vertex operator algebra
since the spaces of weights one of them are abelian and Lie algebra $\fg=D_8$, respectively. However, they have the same central charge $(c=8)$ and conformal weights $\vecr=\{0,\, 1/2,\, 1\}$ and $\vech=\{0,\,1/2,\,4/3\}$. Therefore one can show that the MLDE associated with both VOAs is same, that is given by 
\begin{equation}\label{cc=8.EQ2}
f'''-\frac{1}{2}E_2\,f''+\Bigl(\frac{1}{2}E_2'-\frac{1}{6}E_4\Bigl)f'+\frac{1}{27}E_6\,f\=0\,.
\end{equation}
Moreover,  one can find explicit solutions of~\eqref{cc=8.EQ2} as follows. 

The main idea of finding a~solution of the MLDE~\eqref{cc=8.EQ2} is as follows. First, we suppose by Theorem~\ref{T-DLN} 
that $f_i$ is a modular function and level six since the set of indicial roots is $\left\{-1/3,\,1/6,\,2/3\right\}$. 
Therefore, $\eta^8f_i$ $(1\leq i\leq3)$ would be a modular form of weight 4 on $\Gamma(2)$ and $[\Gamma(2)\,:\,\mathrm{SL}_2(\mathbb{Z})]=6$,\,$\dim_\mathbb{C}M_k(\Gamma(2))=1+[k/2]\,\text{for}\,k\in2\mathbb{Z}_{\ge0}$. 
It is well known that a basis of the space of modular forms of weight 4 on $\Gamma(2)$ is given by 
\begin{equation}\label{q-ser2}
\begin{aligned}
&\eta(q)^8x(q)^2\=1+48q+624q^2+1344q^3+5232q^4+6048 q^5+\cdots\,,\\
&\eta(q)^8x(q)y(q)\=q^{1/2}(1 + 28 q + 126 q^2 + 344 q^3 + 757 q^4+\cdots)\,,\\
&\eta(q)^8y(q)^2\=q+8q^2+28q^3+64q^4+126q^5+\cdots\,,
\end{aligned}
\end{equation}
where $x(q)$ and $y(q)$ are given by 
\begin{equation*}
x(q)\=\frac{2E_2(q^2)-E_2(q)}{\eta(q)^4}\quad\text{and}\quad y(q)\=\frac{\eta(q^2)^4}{\eta(q)^4}\,,
\end{equation*}
respectively, and they have relations
\begin{equation*}
x'(q)\=-\frac{1}{6}\eta(q)^4\bigl(x(q)^2-192y(q)^2\bigr)\quad\text{and}\quad y'(q)\=\frac{1}{3}\eta(q)^4x(q)y(q)\,.
\end{equation*}

Secondly, since $\eta^8f_i\in M_{4}(\Gamma(2))$ for each~$i$ is a linear combination of the functions of~\eqref{q-ser2}, 
comparing $q$-expansions of those two (first three coefficients are enough since $\dim_\C M_{4}(\Gamma(2))=3$), 
for example, one can conclude that $f_2=xy$. In fact, the function $\eta^8f_2$ and $\eta^8xy$ satisfy the same equations 
which are equivalent to~\eqref{thirdMLDE.REC}. Since enough number of Fourier coefficients of $\eta^8f_i$ $(1\leq i\leq3)$ 
and $\eta^8xy$, for instance, are determined by using the recursive relation~\eqref{thirdMLDE.REC} of $f_i$ 
(Fourier expansion of $\eta$ is obvious) and by the definition of $x$ and $y$, respectively, one can verify $\eta^8f_2=\eta^8xy$. 
Then the exactly same discussions as above give~\eqref{q-ser2}.
Summarizing one has 
\begin{equation}
f_1\=x^2\,,\quad f_2\=xy\,,\quad f_3\=y^2\,.
\end{equation}

The conformal weights are 0, $1/2$ and 1 and have the integral difference~1 (between~0 and~1). 
Hence $f_1+sf_3$ for any complex number $s$ is a solution with the indicial root~0. Therefore, for instance,  
$\{f_1+64f_3,\,f_2,\,f_3\}$ is another basis of the space of solutions of~\eqref{cc=8.EQ2}. The number 64 is chosen as 
$f_1+64f_3$ coincides with $\ch_{V_{D_8}}(q)=(\theta_3(q)^{8}+\theta_0(q)^8)/2\eta(q)^8$
Since the dimensions of the weight one spaces of the lattice $\dim V_{D_8}$ and the affine VOA $\AL{D}{8}{1}$ are~120,
the characters of them are same and given by
\begin{equation*}
f_1+64f_3\=x(q)^2+64y(q)^2\=q^{-1/3}(1 + 120 q + 2076 q^2 + 17344 q^3+\cdots)
\end{equation*}
whose Fourier coefficients are non-negative integer since those of~$x(q)$ and~$y(q)$ are non-negative integers. 

The $q$-series $\eta(q)^8f_i$ and $\eta(q)^8$ are modular forms of weight 4 on $\G(2)$ and $\G(3)$, respectively. 
Therefore the functions $f_i$ are modular functions on $\G(6)$ (which does not contradict to Theorem~\ref{T-DLN}). 
Moreover, these solutions are rewritten using theta series as
\begin{equation*}
f_1+64f_3\=\frac{\theta_3(q)^{8}+\theta_0(q)^{8}}{2\eta(q)^8}\,,\quad
f_2\=\frac{\theta_3(q)^{8}-\theta_0(q)^{8}}{2\cdot16\eta(q)^8}\,,\quad
f_3\=\frac{\theta_2(q)^{8}}{2\cdot 128\eta(q)^8}\,
\end{equation*}

Since the $q$-series of indices $-1/3$ and $2/3$ ($f_1+64f_3$, $f_3$, respectively) have an integral difference,
one can find another different basis of the space of solutions of~\eqref{cc=8.EQ2} as 
\begin{align*}
&\frac{E_4(q^2)}{\eta(q)^{8}}\=f_1-48f_3\=q^{-1/3}(1 + 8 q + 284 q^2 + 2112 q^3+\cdots)\,,\\
&x(q)y(q)\=f_2\=q^{1/6}(1 + 36 q + 394 q^2 + 2776 q^3+\cdots)\,,\\ 
&y(q)^2\=f_3\=q^{2/3}(1 + 16 q + 136 q^2 + 832 q^3+\cdots)\,.
\end{align*}
Since the dimensions of lowest weights spaces of simple modules over the lattice VOA $V_{\sqrt{2}\,E_8}$  with conformal weights $0$, $1/2$ and 1 are 1, 16 and 128, respectively, its characters are $f_1-128 f_3$, $16f_2$ and $128f_3$, respectively 
(see \cite[\S~7]{CS} for this fact).

\mn
\DEF
Let $V$~and $W$~be vertex operator algebras.  If $V$~and $W$~are pseudo-isomorphic, we say $V$ and $W$ are in the same
\textit{pseudo-isomorphic class}.

\begin{theorem}\label{T.c=8.1/2}
Let $V$ be a vertex operator algebra with the central charge $c=8$. Suppose the second coefficient of the character 
of $V$ is not $248$. Then  the conformal weights of $V$ is $0$, $1/2$ and $1$, and the space of characters coincides with 
the space of solutions of a~modular linear differential equation of third order. Then $V$~is pseudo-isomorphic to 
the lattice vertex operator algebra $V_{\sqrt{2}\,E_8}$~associated with the integral even lattice~$\sqrt{2}\,E_8$. 
The corresponding MLDE is
\begin{equation*}\label{thirdMLDE5}
f'''-\frac{1}{2}E_2f''+\left(\frac{1}{2}E'_2-\frac{1}{6}E_4\right)f'+\frac{1}{27}E_6\,f\=0\,,
\end{equation*}
and the set of characters of\, $V$~is given by
\begin{equation*}
\ch_V(\tau)\=\frac{E_4(q^2)}{\eta(q)^8}\,,\quad\ch_{1/2}(\tau)\=\frac{\theta_3(q)^{8}-\theta_0(q)^{8}}{2\eta(q)^8}\,,\quad\ch_1(\tau)=4\frac{\theta_2(q)^{8}}{\eta(q)^8}\,.
\end{equation*}
\end{theorem}

\mn
\Rems
(a) Let $V$ be a VOA and $\sigma$ is an automorphism of~$V$. Then it is well known that the fixed point set 
$V^\sigma$ is a vertex operator algebra (which is often called the \textbf{orbifold} of $V$ by $\sigma$). 
Since the structure of orbifold is naturally induced from that of $V$, the central charges of $V$ and $V^\sigma$ are same. 
It is known that the orbifold VOA $V_{\sqrt{2}\,E_8}^+$ has $2^{10}$ simple modules. However, the set of conformal weights
is ${0,\,1,\,1/2}$, and the characters with same conformal weight are equal. The characters of $V_{\sqrt{2}\,E_8}^+$ are given by 
\begin{equation*}
\ch_V(\tau)\=\frac{32E_4(q^2)-\theta_2(q)^8}{32\eta(q)^8}\,,\quad\ch_{1/2}(\tau)=\frac{\theta_3(q)^{8}-\theta_0(q)^{8}}{\eta(q)^8},\,\quad\ch_1(\tau)=8\frac{\theta_2(q)^{8}}{\eta(q)^8}\,.
\end{equation*}
Since
\begin{align*}
&\frac{32E_4(q^2)-\theta_2(q)^8}{32\eta(q)^8}\=q^{-1/3}(1 + 156 q^2 + 1024 q^3 + 6790 q^4 + 32768 q^5+\cdots)\,,\\
&\frac{\theta_3(q)^8-\theta_0(q)^8}{\eta(q)^8}\=32q^{1/6}(1 + 36 q + 394 q^2 + 2776 q^3 + 15155 q^4 +\cdots)\,,\\
&8\frac{\theta_2(q)^8}{2\eta(q)^8}\=2^{11}q^{2/3}(1 + 16 q + 136 q^2 + 832 q^3 + 4132 q^4 +\cdots)
\end{align*}
and then $(-1/3)+(1/6)+(2/3)=1/2=\frac{1}{12}\cdot3\cdot(3-1)$. Therefore $V_{\sqrt{2}\,E_8}^+$ satisfies the non zero Wronskian
condition.Therefore $V_{\sqrt{2}\,E_8}$ is pseudo-isomorphic to $V_{\sqrt{2}\,E_8}^+$, where the automorphism $+$ is defined by 
$\alpha\mapsto-\alpha$ for any $\alpha\in\sqrt{2}E_8$\\
(b) The VOA $V_{\sqrt{2}\,E_8}$ is also pseudo-isomorphic to the affine vertex operator algebra of type $D_{8}$ and level~$1$, whose
characters are given by
\begin{equation*}
\ch_V(\tau)\=\frac{\theta_3(q)^{8}+\theta_0(q)^{8}}{2\eta(q)^8}\,,\quad\ch_{1/2}(\tau)\=\frac{\theta_3(q)^{8}-\theta_0(q)^{8}}{2\eta(q)^8}\,,\quad\ch_{1}(\tau)\=\frac{\theta_2(q)^{8}}{2\eta(q)^8}\,.
\end{equation*} 
which are also written $\ch_{V}(\tau)=f_1-128f_3$, $\ch_{1/2}(\tau)=16f_2$, $\ch_{1}(\tau)=128f_3$ as before.

\mn
\Rems 
(a) We have more two cases, i.e., $h=-1/2, \,2$ and $h=-1/3$.
The former case gives the set of characters of  $V_{D_{20}}$ 
\begin{align*}
&\frac{\theta_3(q)^{20}+\theta_0(q)^{20}}{2\eta(q)^{20}}\=q^{-5/6}(1 + 780 q + 92990 q^2 + 4235960 q^3+\cdots)\\
&\frac{\theta_3(q)^{20}-\theta_0(q)^{20}}{2\eta(q)^{20}}\=40q^{-1/3}\Big(1+248 q + \dfrac{86156}{5} q^2+\cdots\Big)\,,\\
&\frac{\theta_2(q)^{20}}{2\eta(q)^{20}}\=2^{19}q^{5/3}(1 + 40 q + 820 q^2 +\cdots)
\end{align*}
satisfy the MLDE
\begin{equation*}
f'''-\dfrac{1}{2}E_2f''+\Big(\dfrac{1}{2}E_2'-\dfrac{5}{3}E_4\Big)f'-\dfrac{25}{54}E_6 f\=0\,.
\end{equation*}

For $h=-1/3$, there is an ELDE such that  $E_4(q)/\eta(q)^8$ and $E_4(q)^2/\eta(q)^{16}$ are solutions, i.e.
\begin{equation*}
f'''-\dfrac{1}{2}E_2f''+\Big(\dfrac{1}{2}E_2'-\dfrac{23}{18}E_4\Big)f'-\dfrac{1}{3}E_6 f\=0\,.
\end{equation*}
The other solution is 
\begin{equation*}
f_3\=q^{3/2}\Big(1 + \frac{10188}{323}q+\frac{18705546}{37145}q^2+\frac{185597486664}{33393355}q^3
+\frac{1912837788531}{39856585}q^4+\cdots\Big)\,.
\end{equation*}\\
(b) It is proved without any difficulty that the spaces of characters of the lattice VOA $V_{\sqrt{2}\,E_8}$  and the affine VOA 
of type $D_{8}$ with level~1 $(\AL{D}{8}{1})$ are uniquely determined by the condition that 
$a_1=120=\dim\frakg$ where $\frakg=D_8$.\\
(c)~Since conformal weights~0 and~1 have an~integral difference,  the MLDE~\eqref{cc=8.EQ2} might have a~logarithmic solution, 
but it does not happen since the $V_{\sqrt{2}E_8}$ lattice VOA is rational (see~\cite{Z}). Moreover, solutions with the conformal weight~0 
are not unique since we can add any constant multiple of $\ch_3$ to $\ch_V$ as we discussed in the previous paragraphs,
where $ch_*$ is the character of simple $V$-modules with conformal weight $*$ and $\ch_V$ is the character of the vacuum module.

\mn
\Rem
The MLDEs associated with the lattice VOA $V_{\sqrt{2}\,E_8}$ and the affine VOA $L_{D_8,1}$ are equal. 
However $V_{\sqrt{2}\,E_8}$ and $L_{D_8,1}$ are not isomorphic since the weight one spaces of $V_{\sqrt{2}\,E_8}$ 
and $L_{D_8,1}$ are abelian and $\fg=D_8$ by 0th product, respectively. Therefore, we do not know if $V$ is isomorphic to 
either $V_{\sqrt{2}\,E_8}$ or $L_{D_8,1}$. In other words, this give rises to an example such that the space of the characters 
does not determine VOAs uniquely.

\subsection{Exceptional case}\label{SS.m=248}
In this subsection we study the case $m=248$. (See the first paragraph of section~\ref{HGF}).
The MLDE associated with $m=248$ is given by
\begin{equation}\label{cc8-248}
f'''-\frac{1}{2}E_2 f''+\Bigl(\frac{1}{2}E_2'+xE_4\Big)f'+yE_6 f\=0\,.
\end{equation}
Since $f_0=q^{-1/3}(1+248q+\cdots)$ is a~solution of~\eqref{cc8-248} (see~\eqref{E.c8E4}), one has $-5 + 18 x + 54 y=0$.
Then the corresponding MLDE is written as 
\begin{equation}\label{cc8-248.1}
f'''-\frac{1}{2}E_2 f''+\Bigl(\frac{1}{2}E_2'+xE_4\Big)f'+\Big(\frac{x}{3}+\frac{5}{54}\Big)E_6 f\=0\,,
\end{equation}
where $x$ is an~arbitrary complex number. 

First, we suppose that $f_1=q^{-1/3}(1+248q+m_2q^2+\cdots)$ is a~solution of~\eqref{cc8-248}.
Then we have $(5+3 x) (m_2-4124)=0$. If $5+3x=0$, we have $(c;h_1,h_2)=(8,1/2,1)$ since the indicial equation is 
$t^3-t^2/2-t/6+1/27=(t+1/3)(t-1/6)(t-2/3)=0$. Therefore, eq.~\eqref{cc8-248.1} has a basis of the space of solution
$f_1=E_4(q)/\eta(q)^8$, $f_2=(\theta_3(q)^8-\theta_0(q)^8)/2\eta(q)^8$ and $f_3=\theta_2(q)^8/2\eta(q)^8$,  
which are same as being stated in Theorem~\ref{T.c=8.1/2}.

Secondly, we suppose that $m_2=4124$ and $f_1=q^{-1/3}(1+248q+4124q^2+m_3q^3+\cdots)$.
Then it follows that $(31+6 x)(m_3-34752)=0$. If $31+6 x=0$, the indicial equation is $( 3t-5) (3t+1) (6t+5)=0$.
The set of solutions of the corresponding MLDE is $(\theta_3(q)^{20}+\theta_0(q)^{20})/2\eta(q)^{20}$, 
$(\theta_3(q)^{20}-\theta_0(q)^{20})/80\eta(q)^{20}$, $\theta_2(q)^{20}/2\eta(q)^{20}$. 
In particular, one can verify that the second solution has non-integral coefficients as
\begin{equation*}
\frac{\theta_3(q)^{20}-\theta_0(q)^{20}}{80\eta(q)^{20}}\=\frac{1}{q^{1/3}}+248 q^{2/3}+\frac{86156 }{5}q^{5/3}+\cdots\,.
\end{equation*}

Let us suppose that $m_3-34752=0$ and 
\begin{equation*}
f_1\=q^{-1/3}(1+248q+4124q^2+34752q^3+m_4q^4+\cdots)
\end{equation*}
is a~solution of MLDE~\eqref{cc8-248}. Then we have $(32+3 x) (m_4-213126)=0$.

If $m=248$, the $q$-series $f_0:=q^{-1/3}(1+248q+\cdots)$ is a solution of 
\begin{equation*}
f'''-\frac{1}{2}E_2 f''+\Bigl(\frac{1}{2}E_2'+xE_4\Big)f'+yE_6 f\=0
\end{equation*}
so that 
\begin{equation*}\label{E.248}
f'''-\frac{1}{2}E_2 f''+\Bigl(\frac{1}{2}E_2'+xE_4\Big)f'+\Big(\frac{x}{3}+\frac{5}{54}\Big)E_6 f\=0\,.
\end{equation*}
The the set of incidial roots of~\eqref{E.248} is $\{-1/3, \,5\pm\sqrt{3(-48 x-5)})/12\}$.
Then Ramanujan's relations yield that~\eqref{E.248} has the solution $E_4(q)/\eta(q)^8$ for any $x$.
Suppose that $f_0\neq E_4(q)/\eta(q)^8$. Then there exists a pair of  indicial roots, which has an integral difference.
More precisely, there exists $n\geq2$ such that there exists $x$ such that 
\begin{equation*}
(5-\sqrt{3(-48 x-5)})/12=n-1/3\quad\text{or}\quad (5+\sqrt{3(-48 x-5)})/12=n-1/3\,.
\end{equation*}
Therefore we have $x=-(6 n^2-9 n+4)/6$. Then indicial roots of eq.~\eqref{E.248} is $\{(7-6n)/6,\, -1/3,\, (3n-1)/3\}$. 
Since $n\geq2$, the effective central chargeは$\tilde{c}=4 (6n-7)$ and  the set of conformal weights is $\{0,\,(2n-3)/2,\,(4n-3)/2\}$.
When $n=1$ the set of solutions gives the set of characters of $L_{D_8,1}$ with the vacuum character
$E_4(q)/\eta(q)^8$. When $n=2$ we have the set of characters of $L_{D_{20},1}$の $(h=1/2)$.

Now let us suppose that $n\geq3$. Then $q$-expansion of $f_3$ whose index is $(7-6n)/6$ is 
\begin{equation*}
f_3(q)\=q^{7/6-n}\Big(1-\frac{4(6 n-7)(40 n^2-30 n+17)}{(2 n-5) (4 n-5)}q+\cdots\Big)\,,
\end{equation*}
Therefore the second coefficients is non-negative integer if and only if $n\leq 7/6$ or $5/4<n<5/2$.
Hence we have either $n=1$ or~2. 

\section{Vertex operator algebras with central charge $c=16$}\label{S.cc=16} %Section 6
We now discuss the case that $h=-1,\,3/4$. The set of indicial roots of~\eqref{third.EQ.h-1.1} is $\{-5/3, -2/3, 17/6\}$. 

The same method in the previous Section and its Subsections work for $c=16$ to classify the spaces of characters of VOAs 
with central charge~$c=16$. 

Suppose that~\eqref{E.Third.EQ1} with $c=16$ has a solution of the form $f_2=q^\lambda\bigl(\sum_{j=0}^{\infty}b_jq^{j}\bigr)$ 
(with $b_0=1$ and~$\lambda=h-2/3$). Then one can see that the second coefficient is 
\begin{equation*}
b_1\=8(60h^3-277 h^2+ 437 h-186)/(h+1)(4h-3)
\end{equation*}
if  $h\neq-1,\,3/4$. For $h=-1$ we have
\begin{equation}\label{third.EQ.h-1.1}
f'''-\frac{1}{2}E_2f''+\left(\frac{1}{2}E'_2-\frac{2}{11}\right)f'-\frac{85}{27}E_6f\=0\,,
\end{equation}
For $h=3/4$, we have
\begin{equation}\label{third.EQ.h3/4.1}
f'''-\frac{1}{2}E_2f''+\left(\frac{1}{2}E'_2-\frac{11}{16}\right)f'-\frac{13}{216}E_6f\=0\,,
\end{equation}
The set of solutions of~\eqref{third.EQ.h-1.1} is 
\begin{align*}
f_1&\=\dfrac{E_4(q)^2}{\eta(q)^{16}}\=q^{-2/3}(1 + 496 q + 69752 q^2 + 2115008 q^3 +\cdots)\,,\\
f_2&\=q^{-5/3}\Big(1 +\dfrac{1254592}{1617}q -\dfrac{800544692}{1617}q^2 -\cdots\Big)+\dfrac{7680}{7}(\log q)\dfrac{E_4(q)^2}{\eta(q)^{16}}\,,\\
f_3&\=q^{17/6}\Big(1 +\dfrac{2788}{33}q+\dfrac{17400758}{4719}q^2 +\dfrac{6701005192}{61347}q^3 +\cdots\Big)\,.
\end{align*} 
The set of solutions of~\eqref{third.EQ.h3/4.1} is $f_1=E_4(q)^2/\eta(q)^{16}$, 
\begin{align*}
f_2&\=q^{13/12}\Big(1+\frac{416}{11}q+\frac{7709}{11}q^2+\frac{1799980}{209}q^3+\cdots\Big)\,,\\
f_3&\=q^{1/12}\Big(1 +\dfrac{8464}{49}q +\dfrac{33390383}{5929}q^2
+\dfrac{562460144}{5929}q^3 +\cdots\Big)+\dfrac{180}{7}(\log q)f_2\,,
\end{align*}
Therefore there is no VOA whose characters are as above.

Since it is immediate that if $y\,b_1$ is an~integer for a~positive integer~$y$, 
then $n= (b_1-120h+584)y=120y(37h-27)/(h+1)(4h-3)$ is also an integer, and this equation is rewritten as 
$4n h^2+(n-4440y)h-3(n-1080y)=0$ for~$n\in\Z$.  This Diophanus equation in indeterminate $h$ has a rational solution 
if and only if $(n, d, y)\in\Z^3$ satisfies 
\begin{equation*}
(49n-30360y-7d)(49n-30360y+7d)\=-44236800y^2\,.
\end{equation*} 

We now suppose that $y=16$, which is motivated by the fact that the maximal dimension of the lowest weight spaces of modules of the orbifold VOA $V_{\Lambda_{16}}^+$ of the Barnes--Wall lattice $\Lambda_{16}$ is~16. 
It then follows that the value of h is one of an element of the set
\begin{align*}
\{&-4809933/8\,, -719755/4\,, -2130037/16\,, -74613\,, -21950\,, -32197/2\,, -45669/4\,,\\
& -8784\,, -4395\,, -2200\,, -16429/8\,, -1468\,, -882\,, -552\,, 0\,, 280\,, 299\,, 376\,, 456\,, 504\,, \\
&600\,,680\,, 728\,, 792\,, 852\,, 1080\,, 1450\,, 1496\,, 1744\,, 2184\,, 5247/2\,, 4380\,, 29225/4\,,\\
&904031/2\,, 57788591/16\,\}\,.
\end{align*}

Now, one can see that $16b_1\,(=yb_1)$ is a non-negative integer if and only if $h$ is an element in the set 
\begin{equation*}
\Bigl\{\,-\frac{11}{12},\,-\frac{17}{20},\, -\frac{3}{4},\, -\frac{7}{12},\, -\frac{1}{2},\, -\frac{3}{8},\, -\frac{1}{4},\, 0,\, 
\frac{1}{4},\, \frac{1}{3},\, \frac{1}{2},\, \frac{3}{5},\, \frac{2}{3},\, \frac{7}{8},\, 1,\, \frac{3}{2},\, 2,\, \frac{11}{4},\, 3\,\Bigr\}\,. 
\end{equation*}
It hence follows by Lemma~\ref{SYM.INV} that the values of $h$ give a subset
\begin{equation*}
\left\{\,-\frac{1}{2},\,-\frac{1}{4},\,\frac{1}{2},\,1,\,\frac{3}{2},\,2,\,\frac{11}{4},\,3\,\right\}\,.
\end{equation*}

Finally, since the denominator of irreducible fractions of the fourth coefficient\footnote{This coefficient given described in Appendix.} 
of $q$-series $f_2$ has a factor $1+4h$, one must have $h\neq-1/4$. Hence $h\neq\mu(-1/4)=11/4$ for $\mu\in G$. 
Therefore $h$ must be an element of the set $\{1,\,2,\,3\,\}$. We now study the MLDE for each~$h$ case by case basis.

For $h=-1/4$ we have
 \begin{equation}\label{third.EQ.h3/4}
f'''-\frac{1}{2}E_2f''+\left(\frac{1}{2}E'_2-\frac{43}{16}\right)f'-\frac{275}{216}E_6f\=0\,,
\end{equation}
whose the set of incidial roots is $\{-11/12, -2/3, 25/12\}$. The solution of the indicial root $-2/3$ is
$E_4(q)^2/\eta(q)^{16}$. Suppose that $f=q^{-11/12}(a_0+a_1 q+\cdots)$ is a solution of~\eqref{third.EQ.h3/4}.
Then one can have $f=q^{-11/12}(1+836q+\ell q^2+\cdots)$ and $7\ell-1034649=0$ and $17466372 - 56\ell=0$, which
is a contradiction.

\bn
\textbf{(a)} $h=1$. The set of indicial roots of~\eqref{third.EQ.h-1.1} is $\{-5/3, -2/3, 17/6\}$. The set of conformal weights is $\vh=\{0,\,1,\,3/2\}$ which coincides with 
that of the \textit{Barnes--Wall} lattice VOA, and the set of indices is $\vr=\{-2/3, \,1/3, \,5/6\}$. 
Hence the corresponding MLDE is given by 
\begin{equation}\label{cc=16.a}
f'''-\frac{1}{2}E_2f''+\left(\frac{1}{2}E'_2-\frac{1}{2}E_4\right)f'+\frac{5}{27}E_6\,f\=0\,.
\end{equation}
The first several Fourier coefficients of $\eta(q)^{16}f_2(q)$ and $\eta(q)^{16}f_3(q)$ being computed  
by the recursive relations~\eqref{thirdMLDE.REC} are
\begin{equation*}
\begin{aligned}
&\eta(q)^{16}f_2(q)\=q + 120 q^2 + 2060 q^3 + 15424 q^4 + 73518 q^5 + 263584 q^6+\cdots\,,\\
&\eta(q)^{16}f_3(q)\=q^{3/2}\left(1+36q+378 q^2+2200q^3+8955q^4+28836 q^5+\cdots\right)\,.
\end{aligned}
\end{equation*}
Moreover, by Theorem~\ref{T-DLN} the functions $\eta(q)^{16}f_i(q)$ $(i=2,\,3)$ are supposed to be modular forms 
of weight 8 on~$\G(2)$, which can be also seen since $\eta(q)^{16}f_i(q)$ are invariant under the action $\tau\mapsto\tau+2$.

It is known that the space of modular forms of weight~$8$ on~$\Gamma(2)$ without constant terms of their Fourier expansions is 
a subspace of the union of $q\cdot\mathbb{C}[[q]]$ and $q^{1/2}\cdot\mathbb{C}[[q]]$,
and are linearly generated by 
\begin{align}\label{basis_16-1}
\begin{split}
&\eta(q)^{16}x(q)^2y(q)^2\\
&\qquad\=q + 56 q^2 + 1036 q^3 + 7744 q^4 + 36654 q^5 + 131744 q^6+\cdots\,,
\end{split}\\
\begin{split}\label{basis_16-2}
&\eta(q)^{16}y(q)^4\\
&\qquad\=q^2 + 16 q^3 + 120 q^4 + 576 q^5 + 2060 q^6 + 6048 q^7+\cdots\,,
\end{split}\\
\begin{split}\label{basis_16-3}
&\eta(q)^{16}x(q)y(q)^3\\
&\qquad\=q^{3/2}\left(1 + 36 q + 378 q^2 + 2200 q^3 + 8955 q^4 + 28836 q^5 +\cdots\right)\,,
\end{split}\\
\begin{split}
&\eta(q)^{16}x(q)^3y(q)\\
&\qquad\=q^{1/2}\left(1 + 76 q + 2094 q^2 + 25208 q^3 + 138757 q^4 + 574212 q^5+\cdots\right)\,,
\end{split} \label{basis_16-4}
\end{align}
respectively. By comparing the first two coefficients of $\eta(q)^{16}f_2(q)$ and $\eta(q)^{16}f_3(q)$ with a linear combination of~\eqref{basis_16-1} and~\eqref{basis_16-2}, respectively, one can have
\begin{equation*}
f_2\=y^2\bigl(x^2+64y^2\bigr)\quad \text{and}\quad f_3\=xy^3\,.
\end{equation*}

Since indices $-2/3$ and $1/3$ have an integral difference, the space of solutions with the index $-2/3$ is linearly generated by 
$E_4^2/\eta^{16}$ and $f_2$. Let $\Theta_{\Lambda_{16}}(q)$ be the theta series of the Barnes--Wall lattice $\Lambda_{16}$ 
(with rank 16), which is a 16-dimensional positive definite even integral lattice of discriminant~28, more precisely, which does not have vectors with norm~2. (See \cite[Section~4.10]{CS} for more details.)
\begin{equation*}
\begin{split}
\Theta_{\Lambda_{16}}(q)&\=\frac{1}{2}\bigl(\theta_2(q^2)^{16}+\theta_3(q^2)^{16}
+\theta_0(q^2)^{16}+30\theta_2(q^2)^8\theta_3(q^2)^8\bigr)\\
&\=1+4320 q^2+61440q^3+\cdots\,,
\end{split}
\end{equation*} 
which is nothing but eq.~(132) in~\cite{CS}. Since $\Theta_{\Lambda_{16}}(q)=E_4(q)^2-480\eta(q)^{16}f_2(q)$ 
(of course, the levels of  $E_4(q)$ and $\eta(q)^{16}f_2(q)$ are smaller than that of $\Theta_{\Lambda_{16}}(q)$), 
the characters satisfying~\eqref{cc=16.a} are given by 
\begin{align*}
\begin{split}
&\ch_V(\tau)\=\frac{\Theta_{\Lambda_{16}}(q)}{\eta(q)^{16}}\\
&\qquad\quad\=q^{-2/3}(1 + 16 q + 4472 q^2 + 131648 q^3 + 2168860 q^4 + 24647840 q^5+\cdots)\,,
\end{split}\\
\begin{split}
&\ch_1(\tau)\=\frac{\theta_2(q)^8(\theta_3(q)^8+\theta_0(q)^8)}{16\eta(q)^{16}}\\
&\qquad\quad\=32q^{1/3}\left(1+136 q+4132 q^2+67712 q^3+770442q^4+6834240 q^5+\cdots\right)\,,
\end{split}\\
\begin{split}
&\ch_{3/2}(\tau)\=\frac{\theta_2(q)^8(\theta_3(q)^8-\theta_0(q)^8)}{16\eta(q)^{16}}\\
&\qquad\quad\;\;\=512q^{5/6}\left(1 + 52 q + 1106 q^2 + 14808 q^3 + 147239 q^4 + 1183780 q^5+\cdots\right)\,.
\end{split}
\end{align*}
The $q$-series $\eta(q)^{16}f_i(q)$ and $\eta(q)^{16}$ are modular forms of weight 8 on $\G(2)$ and $\G(3)$, respectively, 
and hence these characters are modular functions on~$\G(6)$. (There is another way to verify the modularity of the solutions 
(see~\cite[Theorem~I\,]{DLN}. Thus the value which is the smallest positive integer $N$ such that all exponents
$\{-2/3, 1/3, 5/6\}$ multiplied by $N$ are integers gives rise to the rank of congruence groups, which is~6.)

Finally, the  solution $E_4^2/\eta^{16}$ is rewritten as 
\begin{equation*}
\frac{E_4(q)^2}{\eta(q)^{16}}\=\ch_V(\tau)+\ch_1(\tau)\=q^{-2/3}(1 + 496 q^2 + 2115008 q^3+\cdots)\,.
\end{equation*}

\bn
\textbf{(b) $h=2$.} 
The set of  conformal weights is $\vh=\{0,\,1/2,\,2\}$ and the set of indices is $\vecr=\{-2/3,\, -1/6,\,4/3\}$, 
which coincides with the set of conformal weights of the affine VOA of type $D_{16}$ and level one. 
The corresponding MLDE is given by 
\begin{equation*}
f'''-\frac{1}{2}E_2\,f''+\Bigl(\frac{1}{2}E'_2-E_4\Bigr)f'-\frac{4}{27}E_6\,f\=0\,.
\end{equation*}
By~\eqref{MLDE_for_D}  one can find  
\begin{align}
\begin{split}
&\ch_V(\tau)\=\frac{\Theta_{D_{16}}(q)}{\eta(q)^{16}}\\
&\qquad\quad\=q^{-2/3}(1 + 496 q + 36984 q^2 + 1066432 q^3+\cdots)\,,\label{ch_D_16_1}
\end{split}\\
\begin{split}
&\ch_{1/2}(\tau)\=\frac{\theta_3(q)^{16}-\theta_0(q)^{16}}{2\eta(q)^{16}}\\
&\qquad\qquad=32q^{-1/6}(1 + 156 q + 6790 q^2 + 142136 q^3+\cdots)\,,\label{ch_D_16_2}
\end{split}\\
&\ch_2(\tau)\=\frac{\theta_2(q)^{16}}{2\eta(q)^{16}}
\=32768q^{4/3}(1 + 32 q + 528 q^2 + 6016 q^3+\cdots)\,,  \label{ch_D_16_3}
\end{align}
where $\Theta_{D_{16}}(q)=\big(\theta_3(q)^{16}+\theta_0(q)^{16}\big)/2$ is the lattice theta function of the lattice~$D_{16}$. 
In fact, we have 
\begin{align}
&\eta(q)^{16}f_2(q)\=q^2+16q^3+120q^4+576q^5+2060q^6+\cdots\,,\label{q-ser_D_16_2}\\
&\eta(q)^{16}f_3(q)\=q^{1/2}(1+140q+4398q^2+49400q^3+279557q^4+\cdots)\,.\label{q-ser_D_16_3}
\end{align}

The intersection of the space of modular forms of weight $8$ on $\G(2)$ and $q\,\C[[q]]$ is linearly generated by~\eqref{basis_16-1} 
and~\eqref{basis_16-2}. Similarly, the intersection of the space of modular forms of weight $8$ on $\G(2)$ and $q^{1/2}\C[[q]]$ 
is linearly generated by~\eqref{basis_16-3} and~\eqref{basis_16-4}. Hence by comparing at most first two coefficients. 
(The space of modular forms of weight 8 on $\Gamma(2)$ is 5-dimensional. However, the intersection such modular forms 
and $\C[[q]]$ is 3-dimensional and they are represented as $1+O(q)$, $q+O(q^2)$ and $q^2+O(q^3)$, respectively.
By~\eqref{q-ser_D_16_2} and~\eqref{q-ser_D_16_3} with~\eqref{basis_16-2} and eqs.~\eqref{basis_16-3} and \eqref{basis_16-4}, respectively, we have $f_2=y^4$ and $f_3=xy\bigl(x^2+64y^2\bigr)$. Moreover, using relations $2\eta^4x=\theta_3^4+\theta_0^4$, $2^4\eta^4y=\theta_2^4$ and $\theta_3^4-\theta_0^4=\theta_2^4$, one can verify that $\eta^{16}f_2=\theta_2^{16}/2^{16}$,
$\eta^{16}f_3=\bigl(\theta_3^{16}-\theta_0^{16}\bigr)/32$ 
and $\eta^{16}f_1=\bigl(\theta_3^{16}+\theta_0^{16}\bigr)/2-2^{15}\eta^{16}f_2$. Finally, we have characters (up to constant multiples) 
of the affine VOA $D_{16}$ with level 1 obtained in section~\ref{S.Lattice} in another way. The solution $E_4^2/\eta^{16}$ 
of the MLDE is written as 
\begin{equation*}
\frac{E_4(q)^2}{\eta(q)^{16}}\=\ch_V(\tau)+\ch_2(\tau)\=q^{-2/3}(1 + 496 q + 69752 q^2 + 2115008 q^3+\cdots)\,.
\end{equation*}

\mn
\Rem 
Though the conformal weights 0 and 2 have an~integral difference, the characters are explicitly determined 
by~\eqref{ch_D_16_1}--\eqref{ch_D_16_3} since the second coefficient $\ch_V(\tau)$~coincides with~$\dim D_{16}$ 
(see~(b) of Remarks after Theorem~\ref{T.c=8.1/2} in Section~\ref{S.cc=8}).

\bn
\textbf{(c)~$h=3$.} 
The set of conformal weights is $\vh=\{0,\, 3,\, -1/2\}$ and the set of indices is $\vecr=\{-7/6,\,-2/3,\, 7/3\}$. 
To avoid using modular forms of negative weights (the negative conformal weight $-1/2$ gives rise to a modular form whose
$q$-series is $q^{-1/2}+O(q^{1/2})$ which has a pole at $+i\infty$, i.e., it is a \textit{weakly holomorphic} modular form. 
Because we want to use (ordinary) modular forms to have expression characters
by using a~symmetry, we can suppose that $\mathsfit{c}=c-24h$ or $24h-2c-12=28$ and the set of conformal weights is 
$\vsh=\{0,\,1/2,\, 7/2\}$, which coincides with that of the affine VOA of type $D_{28}$ with level~1. 

The corresponding MLDE is given by 
\begin{equation*}
f'''-\frac{1}{2}E_2\,f''+\Bigl(\frac{1}{2}E'_2-\frac{7}{2}E_4\Bigr)f'-\frac{49}{27}E_6\,f\=0\,.
\end{equation*}
By~\eqref{ch-D1}, \eqref{ch-D2} and~\eqref{ch-D3}, one can find that the characters are given by
\begin{align*}
&\ch_V(\tau)\=\frac{\Theta_{D_{28}}(q)}{\eta(q)^{28}}
\=q^{-7/6}(1 + 1540 q + 370426 q^2 + 34025432 q^3+\cdots)\,,\\
&\ch_{1/2}(\tau)\=\frac{\theta_3(q)^{28}-\theta_0(q)^{28}}{2\eta(q)^{28}}
\=8q^{-1/6}(7 + 3472 q + 488264 q^2 + 31582272 q^3+\cdots)\,,\\
&\ch_{7/2}(\tau)\=\frac{\theta_2(q)^{28}}{2\eta(q)^{28}}
\=134217728 q^{3/7}(1 + 56 q + 1596 q^2 + 30912 q^3+\cdots)\,.
\end{align*}
where $\Theta_{D_{28}}(q)=(\theta_3(q)^{28}+\theta_0(q)^{28})/2$ is the lattice theta series of $D_{28}$.
Now, $q$-series 
\begin{equation}\label{D_28.SOLUTION}
\begin{aligned}
&\eta(q)^{28}f_3(q)\=1 + 1512 q + 327656 q^2 + 24189984 q^3 +\cdots\,,\\
&\eta(q)^{28}f_2(q)\=q^{7/2}(1 + 28 q + 378 q^2 + 3304 q^3 + 21231 q^4+\cdots)\,
\end{aligned}
\end{equation}
are supposed to be desired modular forms of weight $14\,(=28/2)$ and level~2.
The intersection of the space of modular forms of weight~$14$ on~$\G(2)$ and~$\C[[q]]$ 
is linearly generated by 
\begin{align}\label{D_28.BASIS1}
\begin{split}
&\eta^{28}x^7\=1 + 168 q + 12264 q^2 + 508704 q^3 +13172712 q^4+ \cdots\,,\\
&\eta^{28}x^5y^2\=q + 128 q^2 + 6868 q^3 + 200704 q^4 +3499926 q^5 +\cdots\,,\\
&\eta^{28}x^3y^4\=q^2 + 88 q^3 + 3072 q^4 + 55584 +q^5597716 q^6 + \cdots\,,\\ 
&\eta^{28}xy^6\=q^3 + 48 q^4 + 876 q^5 + 9344 q^6 +69282 q^7 +\cdots\,
\end{split}
\end{align}
Similarly, an~intersection of the space of modular forms of weight~$14$ on~$\G(2)$ and~$q^{7/2}\C[[q]]$ 
is linearly generated by
\begin{align}\label{D_28.BASIS2}
&\eta(q)^{28}y(q)^7\=q^{7/2}(1 + 28 q + 378 q^2 + 3304 q^3 + 21231 q^4 + \cdots)\,.
\end{align}

Comparing~\eqref{D_28.SOLUTION} with a~linear combination of~\eqref{D_28.BASIS1} 
and~\eqref{D_28.BASIS2}, and using relations $2\eta^4x=\theta_3^4+\theta_0^4$, $2^4\eta^4y=\theta_2^4$, 
$\theta_3^4-\theta_0^4=\theta_2^4$, one can show that 
\begin{equation*}
2^{27}\eta^{28}f_2\=\theta_2^{28}/2\quad\text{and}\quad \eta^{28}f_3
\=\bigl(\theta_3^{28}+\theta_0^{28}\bigr)/2\,.
\end{equation*}
Moreover, one can obtain 
\begin{equation*}
56\eta^{28}f_1\=\bigl(\theta_3^{28}-\theta_0^{28}-\theta_2^{28}\bigr)/2\,.
\end{equation*} 
The~$\theta_i^{28}$ $(i=2,\,3,\,0)$ and the~$\eta^{28}$ are modular forms of weight 14 on the congruence group $\G(2)$ 
of~$\MG$ ($N=4,\,N_1=2$ by~Table~\ref{THETA_LEVELS} in Section~\ref{S.Lattice}), and $\G(6)$ 
($N=4,\,N_2=6$ by~Table~\ref{ETA_LEVELS} in Section~\ref{S.Lattice}), respectively. 
Therefore these characters are modular functions on~$\G(6)$ (see Table~\ref{CH_LEVELS}).

\mn
\Rem 
The conformal weights~$1/2$ and~$7/2$ have an~integral difference~3. However, the characters of $V_L$ $(L=D_{28})$
are given by~\eqref{ch-D1}, \eqref{ch-D2} and~\eqref{ch-D3}. Moreover, $E_4^2/\eta^{16}$ is expressed as 
\begin{equation*}
E_4(q)^2/\eta(q)^{16}=\left(\ch_{1/2}(\tau)-\ch_{7/2}(\tau)\right)/56\,.
\end{equation*}

For $m=271$ the corresponding MLDE is given by
\begin{equation*}
f'''-\frac{1}{2}E_2 f''+\Bigl(\frac{1}{2}E_2'-\frac{81}{196}E_4\Big)f'-\frac{533}{10584}E_6 f\=0\,,
\end{equation*}
which has solutions of the $q$-series 
\begin{align*}
f_1&\=q^{-1/6}\Big(1 + 271 q + 4076 q^2 + 30862 q^3 +\frac{5029533}{29}q^4+\cdots\Big)\,,\\
f_2&\=q^{-13/42}\Big(1 +\frac{1742}{7}q+\frac{188850}{49}q^2+\frac{10279088}{343}q^3+\cdots\Big)\,,\\
f_3&\=q^{41/42}\Big(1 +\frac{205}{14} q+\frac{5289}{49}q^2+\frac{5921425}{9947}q^3+\frac{186843109}{69629}q^4+\cdots\Big)\,.
\end{align*}
Thus there are no solutions of the vacuum character-type.

\begin{theorem}\label{T.cc=16}
Let $V$ be a simple vertex operator algebra of CFT and finite type with the central charge 16 and rational conformal weights.
Suppose that the set of the characters of simple $V$-modules gives a basis of  the space of solutions 
of a monic modular linear differential equation of third order. Then we have $h=1,\,2,\,3$ and $(c, h)=(16, 1)$, 
$(\mathsfit{c}, \mathsfit{h})=(20,1/2)$ and $(28,1/2)$. Let $\Lambda_{16}$ is the Barnes--Wall lattice.

\sn
\textup{(a)} 
For $h=1$ the set of conformal weighs is $\{0,\, 1,\, 3/2\}$ and the space of characters of $V$~is pseudo-isomorphic to 
the Barnes--Wall lattice vertex operator algebra $V_{\Lambda_{16}}$. The set of characters is given by 
\begin{align*}
&\ch_V(\tau)\=\frac{\Theta_{\Lambda_{16}}(q)}{\eta(q)^{16}}\,,\\
&\ch_{1}(\tau)\=\frac{\theta_2(q)^8\bigl(\theta_3(q)^8+\theta_0(q)^8\bigr)}{16\eta(q)^{16}}\,, \\
&\ch_{3/2}(\tau)\=\frac{\theta_2(q)^8\bigl(\theta_3(q)^8-\theta_0(q)^8\bigr)}{16\eta(q)^{16}}\,,
\end{align*}
where 
\begin{equation*}
\Theta_{\Lambda_{16}}(q)\=\frac{1}{2}\bigl(\theta_2(q^2)^{16}+\theta_3(q^2)^{16}+\theta_0(q^2)^{16}+30\theta_2(q^2)^8\theta_3(q^2)^8\bigr)\,.
\end{equation*}
\\
\textup{(b)} 
For $h=2$ the set of conformal weights is $\{0,\,1/2,\,2\}$ and $V$ is weakly pseudo-isomorphic to the affine VOA of type $D_{16}$ 
and level~1. The set of characters is given by 
\begin{equation*}
\ch_V(\tau)\=\frac{\theta_3(q)^{16}+\theta_0(q)^{16}}{2\eta(q)^{16}}\,,\quad\ch_{1/2}(\tau)
\=\frac{\theta_3(q)^{16}-\theta_0(q)^{16}}{2\eta(q)^{16}}\,,\quad\ch_{2}(\tau)\=\frac{\theta_2(q)^{16}}{2\eta(q)^{16}}\,.
\end{equation*}
\\
\textup{(c)} 
For $h=3$ the set of conformal weights is $\{0,\,1/2,\,2/7\}$ and $V$ is pseudo-isomorphic
to the affine vertex operator algebra of type~$D_{28}$ with level $1$ \textup{(}however, the central charge is~$28$ after a transformation
of $S_3$ \textup{)}. The set of characters is given by 
\begin{equation*}
\ch_V(\tau)\=\frac{\theta_3(q)^{28}+\theta_0(q)^{28}}{2\eta(q)^{28}}\,,
\quad\ch_{1/2}(\tau)\=\frac{\theta_3(q)^{28}-\theta_0(q)^{28}}{2\eta(q)^{28}}\,,
\quad\ch_{7/2}(\tau)\=\frac{\theta_2(q)^{28}}{2\eta(q)^{28}}\,.
\end{equation*}
\end{theorem}

\mn
\Rems
(a) The orbifold $V_{\Lambda_{16}}^+$ has the central charge 16, whose set of characters is given by
\begin{align*}
&\ch_V(\tau)\=\frac{\Theta_{\Lambda_{16}}(q)+\theta_0(q^2)^{16}}{2\eta(q)^{16}}\,,\\ 
&\ch_{1}(\tau)\=\frac{\theta_2(q)^8\bigl(\theta_3(q)^8+\theta_0(q)^8\bigr)}{32\eta(q)^{16}}\,,\\
&\ch_{3/2}(\tau)\=\frac{\theta_2(q)^8(\theta_3(q)^8-\theta_0(q)^8)}{32\eta(q)^{16}}\,.
\end{align*}
(b) By Theorem~\ref{T.cc=16}, the leading coefficients of the characters corresponding to indices $h-1/3$ are~1, 
$2^5$, $2^3\cdot7$, $2^8$, $2^{12}$, $2^{15}$ and $2^{27}$, respectively.  Therefore if we set $y=2^{27}\cdot7$, 
which is the least common multiple of~1, $2^5$, $2^8$, $2^{12}$, $2^{15}$ and~$2^{27}$,  we have the same results
stated in Theorem~\ref{T.cc=16}. \\
(c) The Barnes--Wall lattice is uniquely determined (up to a pseudo-isomorphism) by requiring $a_1=16$ $(h=1)$.

\section*{Appendix 1}
The exact epression of the Fourier coefficient mentioned in the footnote in Section~\ref{S.cc=16} (at the paragraph which contains~\eqref{third.EQ.h3/4}) is given by
\begin{align*}
a^{(2)}_4\=&108840475660 - 17142532250 h + 2159171400 h^2 - 190224000 h^3 + 
 8640000 h^4\\
&-\frac{1279954780160}{11 (h+2)}+\frac{10485088911360}{13 (h+3)}-\frac{416697727057920}{323 (h+4)}+\frac{92570400}{209 (4 h-3)}\\
&-\frac{703761520}{221 (4 h-1)}-\frac{100245600}{4 h+1}+\frac{2471182560}{4 h+3}+\frac{2320465920}{h+1}\,,
 \end{align*}
 with~$h\not\in\{-4,\,-3,\,-2,\,-1,\,-3/4,\,-1/4,\,1/4,\,3/4\,\}$.

\section*{Appendix 2}
In this appendix we study VOAs whose central charge is~4. 

Let $c=4$. As the general form of an~MLDE of third order is 
\begin{equation}\label{E.G3}
f'''-\frac{1}{2}E_2 f''+\Bigl(\frac{1}{2}E_2'+xE_4\Big)f'+y E_6 f\=0\,,
\end{equation}
where $x$ and $y$ is complex numbers. Suppose that $f_1=q^{-1/6}(1+m_1 q+m_2q^2+\cdots)$ is a solution of~\eqref{E.G3}.
The we have 
\begin{align*}
\begin{cases}
&1+ 9 x-54 y\=0\,,\\
&252 + 25 m_1+90 (m_1-48) x+108 (m_1-504) y\=0\,.
\end{cases}
\end{align*}
If $m_1\neq124$ we have a unique solution
\begin{equation}\label{E.xy}
x\=-\frac{1344-48m_1}{192 (124-m_1)}\,,\quad
y\=\frac{160m_1+7808}{6912 (124-m_1)}\,.
\end{equation}

If $m_1=124$, we have $1 + 9 x = 54 y$ and $419 + 855 x - 5130 y=0$, which have no solution $(x,y)$.

We now substitute~\eqref{E.xy} into~\eqref{E.G3}. Then we have 
\begin{equation*}
10m_1^2-(666+m_2)m_1 +8 (17m_2-458)\=0
\end{equation*}
and then 
\begin{equation*}
m_1\=\frac{666+m_2\pm\sqrt{D}}{20}\,,\quad D\=(m_2 - 2054)^2 - 3628800
\end{equation*}
Suppose that there exists $d\in\Z$ such that $D=d^2$, i.e.,
\begin{equation*}
(2054+d-m_2)(2054-d-m_2)\=2^8\cdot3^4\cdot5^2\cdot7
\end{equation*}
Then we have $m_2$ must be one of 210 integers below:
\begin{align*}
&\{-905147, -451548, -300349, -224750, -179391, -149152, -127553, \\
&-111354, -98755, -88676, -73558, -62760, -58441, -54662, -48364, \\
&-43326, -41167, -35770, -34259, -31573, -30374, -28216, -26328, \\
&-23901, -23182, -20666, -19588, -18151, -16894, -16140, -14800, \\
&-14202, -13126, -12409, -12185, -10976, -10618, -10117, -9366, -9227, \\
&-8830, -8116, -7492, -7118, -6691, -6454, -6158, -5626, -5272, -4801, \\
&-4566, -4390, -4144, -3776, -3708, -3514, -3305, -3166, -2935, -2863, \\
&-2682, -2476, -2362, -2220, -2203, -1966, -1798, -1576, -1466, -1384, \\
&-1270, -1141, -1101, -1070, -982, -888, -826, -724, -614, -591, -526, \\
&-478, -419, -412, -316, -250, -199, -166, -126, -113, -97, -58, -16, \\
&6, 32, 35, 58, 74, 98, 120, 124, 134, 140, 146, 149, 3959, 3962, \\
&3968, 3974, 3984, 3988, 4010, 4034, 4050, 4073, 4076, 4102, 4124, \\
&4166, 4205, 4221, 4234, 4274, 4307, 4358, 4424, 4520, 4527, 4586, \\
&4634, 4699, 4722, 4832, 4934, 4996, 5090, 5178, 5209, 5249, 5378, \\
&5492, 5574, 5684, 5906, 6074, 6311, 6328, 6470, 6584, 6790, 6971, \\
&7043, 7274, 7413, 7622, 7816, 7884, 8252, 8498, 8674, 8909, 9380, \\
&9734, 10266, 10562, 10799, 11226, 11600, 12224, 12938, 13335, 13474, \\
&14225, 14726, 15084, 16293, 16517, 17234, 18310, 18908, 20248, 21002, \\
&22259, 23696, 24774, 27290, 28009, 30436, 32324, 34482, 35681, 38367, \\
&39878, 45275, 47434, 52472, 58770, 62549, 66868, 77666, 92784, \\
&102863, 115462, 131661, 153260, 183499, 228858, 304457, 455656, 909255\}
\end{align*}
Since $m_1$and $m_2$ are nonnegative integers, we see that $m_1$ must be one of the following 133 integers:
\begin{align*}
&\{1, 10, 16, 24, 28, 31, 40, 46, 52, 55, 56, 64, 66, 137, 138, 139, 140, 141, 142, 143, 144, 145, 146,\\
&\,\,148, 150, 151, 152, 154, 156, 157, 160, 163, 164, 166, 168, 171, 172, 176, 178, 181, 184, 190, 192,\\
&\,\,196, 199, 206, 208, 216, 217, 220, 226, 232, 241, 244, 248, 256, 262, 271, 276, 280, 296, 298, 304, \\
&\,\,316, 325, 346, 352, 360, 376, 388, 406, 416, 424, 451, 460, 472, 496, 514, 541, 556, 568, 616, 640, \\
&\,\,676, 696, 703, 766, 784, 808, 856, 892, 946, 976, 1000, 1081, 1144, 1216, 1256, 1270, 1396, 1432,\\
&\,\, 1576, 1648, 1756, 1816, 2026, 2152, 2296, 2404, 2656, 2728, 2971, 3160, 3376, 3496, 3916, 4456,\\
&\,\, 4672, 5176, 5806, 6184, 6616, 7696, 9208, 10216, 11476, 13096, 15256, 18280, 22816, 30376, \\
&\,\,45496, 90856\}\,.
\end{align*}

Moreover, since incidial equation
\begin{equation*}
\lambda^3-\frac{1}{2}\lambda^2+x\lambda+y\=0
\end{equation*}
 has only rational solutions, $m_1$ is one of six integers $\{16,\,24,\,28,\,156,\,178,\,271\}$. We now study each case carefully.

\begin{table}[htbp]
\begin{center}
\begin{tabular}{c|c|c}
$m_1$	&$m_2$	&	Indices (Incidial roots)\\ \hline 
$16$		&$98$	&	$(-1/6,1/6,1/2)$\\ \hline
$24$		&$124$	&	$(-1/6,7/30,13/30)$\\ \hline
$28$		&$134$	&	$(-1/6,1/3,1/3)$\\ \hline
$156$	&$6790$	&	$(-2/3,-1/6,4/3)$\\ \hline
$178$	&$4634$	&	$(-1/2,-1/6,7/6)$\\ \hline
$271$	&$4076$	&	$(-13/42,-1/6,41/42)$ \\ \hline
\end{tabular}

\caption{$m_1$, $m_2$ and indices}
\end{center}
\label{default}
\end{table}

\mn
\textbf{Case 1}. $m_1=16$. The corresponding MLDE is 
\begin{equation*}
f'''-\frac{1}{2}E_2 f''+\Bigl(\frac{1}{2}E_2'-\frac{1}{36}E_4\Big)f'+\frac{1}{72}E_6 f\=0
\end{equation*}
which has solutions
\begin{align*}
\frac{I_3(q)^2}{\eta(q)^4}&\=q^{-1/6}(1 + 16 q + 98 q^2 + 364 q^3 + 1221 q^4 +\cdots)\,,\\
\frac{I_3(q)\Delta_3(q)}{\eta(q)^4}&\=q^{1/6}(1 + 11 q + 50 q^2 + 188 q^3 + 583 q^4 +\cdots)\,,\\
\frac{\Delta_3(q)^2}{\eta(q)^4}&\=q^{1/2}(1 + 6 q + 27 q^2 + 92 q^3 + 279 q^4 +\cdots)
\end{align*}
where $I_3(q)=1+6\sum_{n=1}^{\infty}\bigl(\sum_{d | n}\left(\frac{d}{3}\right)\bigr)q^n$ and $\Delta_3(q)=\eta(q^3)^3/\eta(q)$ are modular forms of weight~1 on~$\G(3)$. Also $\eta(q)^4$ is a~modular form of weight~2 on~$\G(6)$. Thus these solutions are modular functions on~$\G(6)$.

\mn
\textbf{Case 2.} $m_1=24$. The MLDE is 
\begin{equation*}
f'''-\frac{1}{2}E_2 f''+\Bigl(\frac{1}{2}E_2'-\frac{1}{100}E_4\Big)f'+\frac{91}{5400}E_6 f\=0
\end{equation*}
whose solutions are
\begin{align*}
\begin{split}
&\psi_1(q)^{10}+14\psi_1(q)^5\psi_2(q)^5-\psi_2(q)^{10}\\
&\qquad\qquad\qquad\qquad\=q^{-1/6}(1 + 24 q + 124 q^2 + 500 q^3 + 1625 q^4 + 4752 q^5 +\cdots)\,,
\end{split}
\\
\begin{split}
&\psi_1(q)^3\psi_2(q)^2\big(\psi_1(q)^5+2\psi_2(q)^5\big)\\
&\qquad\qquad\qquad\qquad\=q^{7/30}\big(1+10 q + 44 q^2 + 164 q^3 + 505 q^4 + 1414 q^5 +\cdots\big)\,,
\end{split}
\\
\begin{split}
&\psi_1(q)^2\psi_2(q)^3\Big(\psi_1(q)^5-\frac{1}{2}\psi_2(q)^5\Big)\\
&\qquad\qquad\qquad\qquad\=q^{13/30}\Big(1 + \frac{13}{2}q+ 30 q^2 +\frac{205}{3}q^3+ 314 q^4 + \frac{1713}{2}q^5+\cdots\Big)\,,
\end{split}
\end{align*}
where 
\begin{equation*}
\psi_1(q)=q^{-1/60}\prod_{\substack{n>0\\ n\equiv\pm2\bmod{5}}}(1-q^n)^{-1}\,,\quad
\psi_2(q)=q^{11/60}\prod_{\substack{n>0\\ n\equiv\pm1\bmod{5}}}(1-q^n)^{-1}
\end{equation*} 
are modular functions on~$\G(60)$. According to the proof o fMilas~\cite{AM3}), all Fourier coefficients of
$\psi_1(q)^5-\psi_2(q)^5$ are non-negative. Hence these solutions do not have negative coefficients.

\mn
\textbf{Cses 3.} $m_1=28$. The corresponding MLDE is 
\begin{equation*}
f'''-\frac{1}{2}E_2 f''+\frac{1}{2}E_2'f'+\frac{1}{54}E_6 f\=0
\end{equation*}
has solutions 
\begin{align*}
&\frac{H_2(q)}{\eta(q)^4}\=q^{-1/6}(1 + 28 q + 134 q^2 + 568 q^3 + 1809 q^4 +\cdots)\,,\\
&\frac{\Delta_2(q)}{\eta(q)^4}\=q^{1/3}(1 + 8 q + 36 q^2 + 128 q^3 + 394 q^4 +\cdots)\,,\\
\begin{split}
&\frac{\Delta_2(q)}{\eta(q)^4}\int_{0}^{q}H_2(q_0)\frac{dq_0}{q_0}+\frac{H_2(q)}{\eta(q)^4}\int_{0}^{q}\Delta_2(q_0)\frac{dq_0}{q_0}\\
&\=16(\log q) q^{1/3}(1 + 8 q + 36 q^2 +\cdots)-q^{1/3}\Big(32+ \frac{1664}{3} q+\frac{33856}{15} q^2+\cdots\Big)\,.
\end{split}
\end{align*}
Therefore there is no VOA whose characters are given above. However, This MLDE is  obtained by acting the Serre derivative 
$\sd_k(f)=f'-k(E_2/12)f$ of weight 4 to the Kaneko-Zagier equation (of $c=4$)
\begin{equation*}
f''-\frac{1}{6}E_2f'-\frac{1}{18}E_4f\=0\,,
\end{equation*}
which is the MLDE whose solutions are characters of the affine (lattice) VOA of type $D_4$ with level~1.

\mn
\textbf{Case 4.} $m_1=156$. The MLDE is
\begin{equation*}
f'''-\frac{1}{2}E_2 f''+\Bigl(\frac{1}{2}E_2'-E_4\Big)f'-\frac{4}{27}E_6 f\=0
\end{equation*}
and the indicial equation is $(t-2/3)(t+1/6)(t-4/3)=0$, whose solutions are 
\begin{align*}
\begin{split}
&\frac{H_2(q)\Delta_2(q)\big(H_2(q)^2+64\Delta_2(q)^2\big)}{\eta(q)^{16}}\\
&\qquad\qquad\qquad\=q^{-1/6}(1 + 156 q + 6790 q^2 + 142136 q^3 + 1897233 q^4 +\cdots)
\end{split}
\\
\begin{split}
&\frac{H_2(q)^4+382H_2(q)^2\Delta_{2A}(q)+14304\Delta_2(q)^4}{\eta(q)^{16}}\\
&\qquad\qquad\qquad\=q^{2/3}(1 + 496 q + 22616 q^2 + 606656 q^3 + 9782812 q^4 +\cdots)\,,
\end{split}\\
&\frac{\Delta_2(q)^4}{\eta(q)^{16}}\=q^{4/3}(1 + 32 q + 528 q^2 + 6016 q^3 + 53384 q^4 +\cdots)\,.
\end{align*}
where $\Delta_{2A}(q)=\eta(q)^8\eta(q^2)^8$ is a~cusp form of weight~8 on~$\G_0(2)$.

\mn
\textbf{Case 5.} $m_1=178$. The MLDE is 
\begin{equation*}
f'''-\frac{1}{2}E_2 f''+\Bigl(\frac{1}{2}E_2'-\frac{25}{36}E_4\Big)f'-\frac{7}{72}E_6 f\=0
\end{equation*}
has solutions
\begin{align*}
\begin{split}
&\frac{I_3(q)^2\Delta_3(q)\big(I_3(q)^3+135\Delta_3(q)^3\big)}{\eta(q)^{12}}\\
&\qquad\qquad\qquad\=q^{-1/6}(1 + 178 q + 4634 q^2 + 61924 q^3 + 566277 q^4 +\cdots)\,,
\end{split}
\\
\begin{split}
&\frac{I_3(q)^6+270\Delta_{3A}(q)+5832\Delta_3(q)^6}{\eta(q)^{12}}\\
&\qquad\qquad\qquad\=q^{-1/2}(1 + 318 q + 8514 q^2 + 126862 q^3 + 1269771 q^4 +\cdots)\,,
\end{split}\\
&\frac{I_3(q)\Delta_3(q)^5}{\eta(q)^{12}}\=q^{7/6}(1 + 23 q + 272 q^2 + 2286 q^3 + 15318 q^4 + \cdots)\,,
\end{align*}
where $\Delta_{3A}(q)=\eta(q)^6\eta(q^3)^6$ is a~cusp form of weight~6 on~$\G_0(3)$.

\mn
\textbf{Case 6.} $m_1=271$. The MLDE
\begin{equation*}
f'''-\frac{1}{2}E_2 f''+\Bigl(\frac{1}{2}E_2'-\frac{81}{196}E_4\Big)f'-\frac{533}{10584}E_6 f\=0
\end{equation*}
which has solutions whose $q$-series has the form
\begin{align*}
f_1&\=q^{-1/6}\Big(1 + 271 q + 4076 q^2 + 30862 q^3 +\frac{5029533}{29}q^4+\cdots\Big)\,,\\
f_2&\=q^{-13/42}\Big(1 +\frac{1742}{7}q+\frac{188850}{49}q^2+\frac{10279088}{343}q^3+\cdots\Big)\,,\\
f_3&\=q^{41/42}\Big(1 +\frac{205}{14} q+\frac{5289}{49}q^2+\frac{5921425}{9947}q^3+\frac{186843109}{69629}q^4+\cdots\Big)\,.
\end{align*}
Thus there are no solutions of the vacuum character-type so that there is no corresponding VOA.

\end{document}